\documentclass{siamltex}
\usepackage{amsmath}
\usepackage{amsfonts}
\usepackage{amssymb}
\usepackage{graphicx}
\usepackage{color}
\usepackage{bm}
\usepackage{verbatim}

\allowdisplaybreaks

\newtheorem{remark}[theorem]{Remark}
\newtheorem{assumption}[theorem]{Assumption}

\newcommand{\be}{\begin{equation}}
\newcommand{\ee}{\end{equation}}
\newcommand{\bea}{\begin{eqnarray}}
\newcommand{\eea}{\end{eqnarray}}
\newcommand{\beas}{\begin{eqnarray*}}
\newcommand{\eeas}{\end{eqnarray*}}

\input{tcilatex}

\newcommand{\tepsilon}{\tilde{\epsilon}}

\begin{document}
\title{A Leray regularized ensemble-proper orthogonal decomposition method for parameterized convection-dominated flows}

\author{Max Gunzburger, Traian Iliescu, and Michael Schneier}

%\date{\today}
\maketitle

\begin{abstract}
Partial differential equations (PDEs) are often dependent on input quantities which are inherently uncertain.
To quantify this uncertainty, these PDEs must be solved over a large ensemble of parameters.
Even for a single realization this can a computationally intensive process.
In the case of flows governed by the Navier-Stokes equations, an efficient method has been devised for computing an ensemble of solutions.
To further reduce the computational cost of this method, an ensemble proper orthogonal decomposition (POD) method was recently proposed.

The main contribution of this work is the introduction of POD spatial filtering for ensemble-POD methods.
The POD spatial filter makes possible the construction of the Leray ensemble-POD model, which is a regularized reduced order model for the numerical simulation of convection-dominated flows.
The Leray ensemble-POD model employs the POD spatial filter to smooth (regularize) the convection term in the Navier-Stokes equations and greatly diminishes the numerical inaccuracies produced by the ensemble-POD method in the numerical simulation of convection-dominated flows.
Specifically, for the numerical simulation of a convection-dominated two-dimensional flow between two offset cylinders, we show that the Leray ensemble-POD method yields accurate results, whereas the ensemble-POD is highly inaccurate.

The second contribution of this work is a new numerical discretization of the variable viscosity ensemble algorithm in which the average viscosity is replaced with the maximum viscosity.
It is shown that this new numerical discretization is significantly more stable than those in current use.
Furthermore, error estimates for the novel Leray ensemble-POD algorithm with this new numerical discretization are also proven.
\end{abstract}


\section{Introduction}

The mathematical models used in realistic applications often times rely on input quantities which are subject to a degree of uncertainty. Some of these quantities include the initial conditions, forcing functions, model coefficients, and the boundary conditions. In order to develop robust models the impact of this uncertainty must be quantified.
%Mathematical flow models inherently rely on fundamental input quantities whose actual values are imperfectly measured or sensitive to variations in other parameters. These quantities include the initial conditions, forcing functions, and model coefficients. Robust flow models must identify and properly handle parameter influence whose sensitivity or uncertainty affects model behavior.  
%In this uncertainty quantification setting, one needs to determine many realizations of the outputs of the simulation in order to determine accurate statistical information about those outputs. A similar need occurs in other settings such as inference, optimization, and control and the need to determine solution ensembles arising in many applications, e.g., weather forecasting and turbulence modeling. 

Common approaches for recovering accurate solutions of these models are the Monte Carlo and stochastic collocation methods \cite{gunzburger_webster_zhang_2014}. These algorithms all require the underlying model to be solved over an ensemble of parameters. Depending upon the problem the spatial resolution required for accurate realizations of the model can render these approaches computationally intractable. In particular, realizations for flow models such as the incompressible Navier-Stokes equations (NSE) can take on the order of weeks. In this work, we are interested in computing ensembles of solutions for the NSE with uncertainty present in the initial conditions, viscosities, and body forces. Specifically, for $j=1,\ldots,J$, we have
\begin{equation}\label{eq:NSE}
\left\{\begin{aligned}
u_{t}^j+u^{j}\cdot\nabla u^{j}-\nu_j\Delta u^{j}+\nabla p^{j}  &
=f^{j}(x,t)&\quad\forall x\in\Omega\times(0,T]\\
\nabla\cdot u^{j}  &  =0&\quad\forall x\in\Omega\times(0,T]\\
u^{j}  &  =0&\quad\forall x\in\partial\Omega\times(0,T]\\
u^{j}(x,0)  &  =u^{j,0}(x)&\quad\forall x\in\Omega,
\end{aligned}\right.
\end{equation}
where $\Omega \subset \mathbb{R}^{d}$, $d= 2,3$, is an open regular domain.

Historically, the solution of the NSE for each parameter has been treated as a separate problem. Recently new algorithms have been developed \cite{J15, J17, JL14,LW17,AMJ16,MohebujjamanR17} that allow for simultaneous calculations at each time step. Specifically the focus of these algorithms has been to use the same linear system for each right hand side. Taking advantage of this problem structure, efficient block solvers, such as block CG \cite{FOP95}, block QMR \cite{FM97}, and block GMRES \cite{GS96} can then be utilized. 
%In most settings, in order to guarantee a desired accuracy level in the outputs, a fine spatial resolution is usually required, which renders each realization to be computationally intensive. Traditionally, the simulation for each ensemble member is treated as a separate problem; therefore the focus has been to cut down on the total number of realizations needed. However, in recent works \cite{J15, J17, JL14,AMJ16}, new algorithms were designed that allow for the realizations to be computed simultaneously at each time step. In those papers, the focus has been on creating algorithms which allow for the same linear system to be used for all right-hand sides. Hence, the need to solve a different linear system for each right-hand side is reduced to solving the a single system with many different right-hand sides. This is a well studied problem for which efficient block iterative methods already exist. A few examples of these include block CG \cite{FOP95}, block QMR \cite{FM97}, and block GMRES \cite{GS96}.

To further improve the efficiency of these ensemble algorithms, reduced order models (ROMs)~\cite{hesthaven2015certified,quarteroni2015reduced} were recently utilized~\cite{GJS17,GJS16}.  Specifically, the proper orthogonal decomposition (POD) method was used to extract the dominant (most energetic) modes from a high-resolution numerical simulation, and the NSE were projected onto these POD modes to obtain an ensemble-POD model.  In~\cite{GJS17,GJS16}, it was shown that the ensemble-POD model significantly decreased the computational cost of the standard ensemble methods, without compromising their numerical accuracy.  We note, however, that the numerical investigation of the ensemble-POD model in~\cite{GJS17,GJS16} was restricted to low Reynolds numbers.  It is well known that, for convection-dominated flows, standard ROMs generally yield inaccurate results, usually in the form of spurious numerical oscillations (see, e.g.,~\cite{giere2015supg,xie2017approximate,xie2017data}).
To mitigate these ROM inaccuracies, several numerical stabilization techniques have been proposed over the years (see, e.g.,~\cite{balajewicz2013low,balajewicz2016minimal,benosman2016robust,carlberg2013gnat,kalashnikova2010stability,osth2014need,wang20162d,2017arXiv170900243C,2017arXiv171003569F, weller2009numerical,weller2009robust}). {\ Regularized ROMs (Reg-ROMs)} are recently proposed stabilized ROMs for the numerical simulation of convection-dominated flows, both deterministic~\cite{sabetghadam2012alpha,wells2017evolve,XIE201812} and stochastic~\cite{iliescu2017regularized}.
These Reg-ROMs use {explicit ROM spatial filtering} to regularize (smooth) various ROM terms and thus increase the numerical stability of the resulting ROM.
This idea was first used by the great Jean Leray~\cite{leray1934sur} in the mathematical study of the NSE and later on in, e.g.,~\cite{geurts2003regularization,layton2012approximate} to develop regularized models for the numerical simulation of turbulent flows~\cite{geurts2003regularization,layton2012approximate}.
In the ROM arena, Reg-ROMs were also successfully used in the numerical simulation of convection-dominated flows.
For example, the Reg-ROMs used in the numerical simulation of a 3D flow past a circular cylinder at a Reynolds number $Re=1000$ produced accurate results in which the spurious numerical oscillations of standard ROMs were significantly decreased~\cite{wells2017evolve}. 

In this paper, we put forth ROM spatial filtering and Reg-ROMs as a means to mitigate the numerical inaccuracies that are generally produced by the ensemble-POD method when this is applied to convection-dominated flows.
Specifically, we propose and investigate the Leray ensemble-POD method, which replaces the convective field in the nonlinearity of the standard ensemble-POD method with its spatially filtered version.
For the spatial filter in the Leray ensemble-POD method, we use the POD differential filter~\cite{wells2017evolve,xie2017approximate}. 
In Section~\ref{POD_sec}, we also propose a new numerical discretization of the variable viscosity ensemble algorithm in which the average viscosity is replaced with the maximum viscosity. 
We show that this new numerical discretization is significantly more stable than those in current use.
Furthermore, we prove error estimates for the new Leray ensemble-POD algorithm with this new numerical discretization.
Finally, in Section~\ref{numex}, we test the new Leray ensemble-POD method in the numerical simulation of the two-dimensional flow between offset circles used in \cite{GJS16,GJS17}.
To this end, we compare the new Leray ensemble-POD method with the standard ensemble-POD method and a fine resolution numerical simulation, which is used as a benchmark.

\section{Notation and preliminaries}

We denote by $\|\cdot\|$ and $(\cdot,\cdot)$ the $L^{2}(\Omega)$ norm and inner product, respectively, and by $\|\cdot\|_{L^{p}}$ and $\|\cdot\|_{W_{p}^{k}}$ the $L^{p}(\Omega)$ and Sobolev
$W^{k}_{p}(\Omega)$ norms, 
respectively. $H^{k}(\Omega)=W_{2}^{k}(\Omega)$ with
norm $\|\cdot\|_{k}$. For a function $v(x,t)$ that is well defined on $\Omega \times [0,T]$, %and for $1 \leq m < \infty$ 
we define the norms
$$
|||v|||_{2,s} : = \Big(\int_{0}^{T}\|v(\cdot,t)\|_{s}^{2}dt\Big)^{\frac{1}{2}}
\qquad \text{and} \qquad 
\||v|||_{\infty,s} := \text{ess\,sup}_{[0,T]}\|v(\cdot,t)\|_{s} .
$$
The space $H^{-1}(\Omega)$ denotes the dual space of bounded linear functionals defined on $H^{1}_{0}(\Omega)=\{v\in H^{1}(\Omega)\,:\,v=0 \mbox{ on } \partial\Omega\}$; this space is equipped with the norm
$$
\|f\|_{-1}=\sup_{0\neq v\in X}\frac{(f,v)}{\| \nabla v\| } 
\quad\forall f\in H^{-1}(\Omega).
$$

The solutions spaces $X$ for the velocity and $Q$ for the pressure are respectively defined as
$$
\begin{aligned}
X : =& [H^{1}_{0}(\Omega)]^{d} = \{ v \in [L^{2}(\Omega)]^{d} \,:\, \nabla v \in [L^{2}(\Omega)]^{d \times d} \ \text{and} \  v = 0 \ \text{on} \ \partial \Omega \} \\
Q : =& L^{2}_{0}(\Omega) = \Big\{ q \in L^{2}(\Omega) \,:\, \int_{\Omega} q dx = 0 \Big\}.
\end{aligned}
$$
A weak formulation of (\ref{eq:NSE}) is given as follows: for $j=1, \ldots, J$, find $u^j:(0,T]\rightarrow X$ and $p^j:(0,T]\rightarrow Q$ such that, for almost all $t\in(0,T]$, satisfy 
\begin{equation}\label{wfwf}
\left\{\begin{aligned}
(u_{t}^j,v)+(u^{j}\cdot\nabla u^{j},v)+\nu_j(\nabla u^{j},\nabla v)-(p^{j}
,\nabla\cdot v)  &  =(f^{j},v)&\quad\forall v\in X\\
(\nabla\cdot u^{j},q)  &  =0&\quad\forall q\in Q\\
u^{j}(x,0)&=u^{j,0}(x).&
\end{aligned}\right.
\end{equation}
The subspace of $X$ consisting of weakly divergence-free functions is defined as
$$
V :=\{v\in X \,:\,(\nabla\cdot v,q)=0\,\,\forall q\in Q\} \subset X.
$$

We denote conforming velocity and pressure finite element spaces based on a regular triangulation of $\Omega$ having maximum triangle diameter $h$ by
$
X_{h}\subset X$ {and} $ Q_{h}\subset Q.
$
We assume that the pair of spaces $(X_h,Q_h)$ satisfy the discrete inf-sup (or $LBB_h$) condition required for stability of finite element approximations; we also assume that the finite element spaces satisfy the approximation properties
$$
\begin{aligned}
\inf_{v_h\in X_h}\| v- v_h \|&\leq C h^{s+1}&\forall v\in [H^{s+1}(\Omega)]^d\\
\inf_{v_h\in X_h}\| \nabla ( v- v_h )\|&\leq C h^s&\forall v\in [H^{s+1}(\Omega)]^d\\
\inf_{q_h\in Q_h}\|  q- q_h \|&\leq C h^s&\forall q\in H^{s}(\Omega),
\end{aligned}
$$
where $C$ is a positive constant that is independent of $h$. The Taylor-Hood element pairs ($P^s$-$P^{s-1}$), $s\geq 2$, are one common choice for which the $LBB_h$ stability condition and the approximation estimates hold \cite{GR79, Max89}.

%In this paper, we solve the NSE \eqref{eq:NSE} (to generate snapshots) using a second-order time stepping scheme (e.g., Crank-Nicolson). 

To ensure the uniqueness of the NSE solution and ensure that standard finite element error estimates, we make the following regularity assumptions on the data and true solution:
{
\begin{assumption}\label{assumption:reg}
In \eqref{eq:NSE} we assume that $u^0 \in V$, $f^{j} \in L^{2}(0,T;L^{2}(\Omega))$, $u^{j} \in L^{\infty}(0,T;H^{s+1}(\Omega))\cap H^{1}(0,T;H^{s+1}(\Omega))\cap
H^{2}(0,T;L^{2}(\Omega))$, and $p \in L^{\infty}(0,T; Q \cap H^k(\Omega))$.
\end{assumption}
}

Using the regularity assumptions above and assuming a sufficiently small $\Delta t$, the following error estimate can be proven for the full discretization of~\eqref{wfwf} with Taylor-Hood elements and the Crank-Nicolson time-discretization~\cite{john2016divergence,layton2008numerical}:
\begin{eqnarray}
	\| u(t^{N}) - u_{h}^{N} \|^2
	+ \nu \, \Delta t \, \sum_{n=1}^{M} \| \nabla (u(t^{n}) - u_{h}^{n}) \|^2
	\leq C \, \left( 
	h^{2m} 
	+ \Delta t^4 
	\right), \, 
	\label{eqn:error-estimate-fe}
\end{eqnarray}
where $C$ is independent of $h$ and $\Delta t$.

We define the trilinear form
$$
b(w,u,v) = (w\cdot\nabla u,v) 
\qquad\forall u,v,w\in [H^1(\Omega)]^d
$$
and the explicitly skew-symmetric trilinear form given by 
$$
b^{\ast}(w,u,v):=\frac{1}{2}(w\cdot\nabla u,v)-\frac{1}{2}(w\cdot\nabla v,u)
\qquad\forall u,v,w\in [H^1(\Omega)]^d \, ,
$$
which satisfies the bounds \cite{Layton08}
\begin{gather}
b^{\ast}(w,u,v)\leq C_{b^*} \|  \nabla w\|   \| \nabla u\| (\|  v \|  \| \nabla
v \| )^{1/2}\qquad\forall u, v, w \in X \label{In1}\\
b^{\ast}(w,u,v)\leq C_{b^*} (\| w \| \|  \nabla w\| )^{1/2}  \| \nabla u\|  \| \nabla
v \| \qquad\forall u, v, w \in X .\label{In2}
\end{gather}
We also define the discretely divergence-free space $V_h$ as
$$
V_{h} :=\{v_{h}\in X_{h}\,:\,(\nabla\cdot v_{h},q_{h})=0\,\,\forall
q_{h}\in Q_{h}\}  \subset X.
$$
In most cases, and for the Taylor-Hood element pair in particular, $V_{h} \not\subset V$, i.e., discretely divergence-free functions are not weakly divergence-free.

\begin{definition}\label{def21}
Let $t^{n}=n\Delta t$, $n=0,1,2,\ldots,N$, where $N:=T/\Delta t$, denote a partition of the interval $[0,T]$. For $j=1, \ldots, J$ and $n=0,1,2,\ldots,N$, let $u^{j,n}(x):=u^{j}(x,t^{n})$. Then, the \text{\bf ensemble mean} is defined, for $n=0,1,2,\ldots,N$, by
$$
 <u>^n : =\frac{1}{J}\sum_{j=1}^{J}u^{j,n}.\label{Enmean}
$$
\end{definition}

The full space and time model which we will base our method off of is similar to the one used in \cite{GJW17,2017arXiv170604060G}. For $j = 1, \ldots,J$, given $u^{j,0}_h \in X_h$ and $u^{j,1}_h\in X_h$, for $n=0,1,2,\ldots,N-1$ find $u^{j,n+1}_h\in X_h$ and $p_h^{j,n+1}\in Q_h$ satisfying
\begin{equation*}%\label{EnB}\emph{}
\begin{aligned}
&\Big(\frac{u^{j,n+1}_h - u^{j,n}_h}{\Delta t}, v_h \Big) + b^{\ast}(<u_h>^{n} , u^{j,n+1}_h ,v_h)+ b^{\ast}(u^{j,n}_h - <u_h>^{n} ,u^{j,n}_h, v_h)\\
& + \nu_{max} (\nabla u^{j,n+1}_h, \nabla v_h) + (\nu_j - \nu_{max}) (\nabla u^{j,n}_h, \nabla v_h) \\
&- (p^{j,n+1}_h , \nabla \cdot v_h)  =( f^{j,n+1}, v_h) \quad \forall v_h\in X_h\\
&\qquad\quad (\nabla \cdot u_h^{j,n+1}, q_h )= 0 \qquad \forall q_h\in Q_h.
\end{aligned}
\end{equation*}

The major difference between the two algorithms the use of maximum value of the viscosities $\nu_{max}$ rather than the average $<\nu>$ resulting in a superior stability condition.

\section{Proper Orthogonal Decomposition Ensemble Based Models}
\label{POD_sec}
\subsection{Proper Orthogonal Decomposition}
In this subsection we briefly describe the POD method and apply it to the previously stated ensemble algorithm. A more detailed description of this method can be found in \cite{KLV01}.

Given a positive integer $N_S$, let $0=t_0<t_1< \cdots < t_{N_S} = T$ denote a uniform partition of the time interval $[0,T]$. For $j=1,\ldots,J_S$, we select $J_S$ different initial conditions $u^{j,0}(x)$, viscosities $\nu^j$, and forcing functions $f^j$  denoted by $u_{h,S}^{j,m}(x)\in X_h$, $j=1,\ldots,J_S$, $m=1,\ldots,N_S$, the finite element approximation to \eqref{eq:NSE} evaluated at $t=t_m$, $m=1,\ldots,N_S$. We then define the space spanned by the $J_S(N_S+1)$ discrete snapshots as
\begin{equation*}
X_{h,S}:=\text{span} \{  u_{h,S}^{j,m}(x) \}_{j=1,m=0}^{J_S,N_S} \subset V_h \subset X_h.
\end{equation*}

Denoting by $\vec{u}_S^{j,m}$ the vector of coefficients corresponding to the finite element function $u_{h,S}^{j,m}(x)$, where $K=\dim X_h$, we define the $K\times J_S(N_S+1)$ {\em snapshot matrix} $\mathbb{A}$ as
$$
\mathbb{A} = \big(\vec{u}_S^{1,0},\vec{u}_S^{1,1}, \ldots , \vec{u}_S^{1,N_S}, \vec{u}_S^{2,0},\vec{u}_S^{2,1},  \ldots , \vec{u}_S^{2,N_S}, \ldots , \vec{u}_S^{J_S,0},\vec{u}_S^{J_S,1}, \ldots , \vec{u}_S^{J_S,N_S}\big),
$$
i.e., the columns of $\mathbb{A}$ are the finite element coefficient vectors corresponding to the discrete snapshots. The POD method then seeks a low dimensional basis 
$$
X_R :=\text{span}\{{\varphi}_i\}_{i=1}^R \subset X_{h,S} \subset V_h\subset X_h
$$
which can approximate the snapshot data. This basis can be determined by solving the constrained minimization problem
\begin{equation}\label{Min}
\begin{aligned}
\min  \sum_{k=1}^{J_S} \sum_{l=0}^{N_S}\Big \|  u_{h,s}^{k,l}-\sum_{j=1}^R (u_{h,s}^{k,l}, \varphi_j)\varphi_j\Big \| ^2 \\
\text{subject to } (\varphi_i, \varphi_j)= \delta_{ij}\quad\mbox{for $i,j=1,\ldots,R$},
\end{aligned}
\end{equation}
where $\delta_{ij}$ denotes the Kronecker delta. Defining the correlation matrix  $\mathbb{C} = \mathbb{A}^{T}\mathbb{M}\mathbb{A}$ where $\mathbb{M}$ denotes the finite element mass matrix, this problem can then be solved by considering the eigenvalue problem
\begin{equation*}
\mathbb{C}\vec{a}_{i} = \lambda_{i}\vec{a}_{i}.
\end{equation*}
It can then be shown the POD basis functions will be given by
\begin{equation*}
\vec\varphi_i = \frac{1}{\sqrt{\lambda_i}}\mathbb{A}\vec{a}_{i}, \ \ \ i = 1, \ldots, R.
\end{equation*}
We now define the POD $L^{2}$ projection we will need for the ensuing stability and error analysis. 
\begin{definition}[POD $L^{2}$ projection]
Let $P_{r}: L^{2}(\Omega) \rightarrow X_{R}$ such that
\begin{equation}
(u - P_{r}u,\varphi) = 0 \qquad\forall \varphi\in X_{R}.
\end{equation}
\end{definition} 
Next we give a POD inverse estimate. Let $\mathbb{S}_{R} = (\nabla \varphi_{i}, \nabla \varphi_{j})_{L^{2}}$ be the POD stiffness matrix and let $\||\cdot|\|_{2}$ denote the matrix $2$-norm
\begin{lemma}[POD inverse estimate]
\begin{equation}\label{POD:inveq}
\|\nabla \varphi \| \leq \||\mathbb{S}_{R}|\|_{2}^{\frac{1}{2}}\|\varphi\| \ \ \ \forall \varphi \in X_{R}.
\end{equation}
\end{lemma}

\subsection{Ensemble-POD Algorithm}
Using this POD basis we can now construct the ensemble-POD algorithm. The construction is similar to  the full finite element approximation except we seek a solution in the POD space $X_{R}$ using the basis $\{ \varphi_{i}\}_{i=1}^{R}$. The fully discrete algorithm can be written as:
\begin{equation}\label{En-POD-Weak}
\begin{aligned}
&\big(\frac{u_{R}^{j,n+1}-u_{R}^{j,n}}{\Delta t}, \varphi\big)+b^{\ast}({<u_{R}
>^{n}},u_{R}^{j,n+1},\varphi)+b^{\ast}({u_{R}^{j, n}-<u_{R}
>^{n}},u_{R}^{j,n}
,\varphi)\\&
+\nu_{max}(\nabla u_{R}^{j,n+1},\nabla
\varphi) + (\nu_{j} - \nu_{max})(\nabla u_{R}^{j,n},\nabla
\varphi) =(f^{j,n+1},\varphi)\qquad\forall \varphi\in X_{R}.
\end{aligned}
\end{equation}
We note that because $X_{R} \subset V_{h}$ the POD basis is
discretely divergence-free by construction. Therefore, there is no pressure term present
in (3.2). In recent works constructing a basis for the pressure space in addition to the velocity space has been investigated. The interested reader should consult \cite{NME:NME4772}.

\subsection{Leray Ensemble-POD Algorithm}
To construct the Leray ensemble-POD algorithm we use the ROM differential filter.

\begin{definition}[ROM differential filter]
$\forall v \in X$ let $\overline{v}^{R}$ be the unique element of $X_{R}$ such that
\begin{equation}
\delta^{2}(\nabla \overline{v}^{R},\nabla \varphi) +  (\overline{v}^{R},\varphi) = (v,\varphi) \qquad\forall \varphi\in X_{R}.
\end{equation}
\end{definition}
Here $\delta$ is known as the filtering radius. The differential filter was first developed by Germano~\cite{germano1986differential} for large eddy simulations. It was introduced in the ROM setting in \cite{sabetghadam2012alpha} and expanded further in \cite{wells2017evolve,XIE201812}. 

Incorporating this into the ensemble framework the fully discrete  Leray ensemble-POD algorithm can be written as:
\begin{equation}\label{En-Leray-POD-Weak}
\begin{aligned}
&\big(\frac{u_{R}^{j,n+1}-u_{R}^{j,n}}{\Delta t}, \varphi\big)+b^{\ast}(\overline{<u_{R}
>^{n}},u_{R}^{j,n+1},\varphi)+b^{\ast}(\overline{u_{R}^{j, n}-<u_{R}
>^{n}},u_{R}^{j,n}
,\varphi)\\&
+\nu_{max}(\nabla u_{R}^{j,n+1},\nabla
\varphi) + (\nu_{j} - \nu_{max})(\nabla u_{R}^{j,n},\nabla
\varphi) =(f^{j,n+1},\varphi)\qquad\forall \varphi\in X_{R}.
\end{aligned}
\end{equation}

%The fully discrete algorithm can be written as:
%\begin{equation}\label{En-Leray-POD-Weak}
%\begin{aligned}
%&\big(\frac{u_{R}^{j,n+1}-u_{R}^{j,n}}{\Delta t}, \varphi\big)+b^{\ast}(\overline{<u_{R}
%>^{n}},u_{R}^{j,n+1},\varphi)+b^{\ast}(\overline{u_{R}^{j, n}-<u_{R}
%>^{n}},u_{R}^{j,n}
%,\varphi)\\&
%+\nu_{max}(\nabla u_{R}^{j,n+1},\nabla
%\varphi) + (\nu_{j} - \nu_{max})(\nabla u_{R}^{j,n},\nabla
%\varphi) =(f^{j,n+1},\varphi)\qquad\forall \varphi\in X_{R}.
%\end{aligned}
%\end{equation}

\section{Stability Analysis}
In this section we present a result pertaining to the stability of the Leray ensemble-POD algorithm. A stability bound for the ensemble-POD algorithm for a fixed viscosity was proven in Theorem 4.2 in~\cite{GJS17}, while a stability bound for an ensemble-FE algorithm with variable viscosity was proven in Theorem 2.1 in~\cite{GJW17}. The stability bound proven in this section is less restrictive than the bound proven in ~\cite{GJW17} due to the use of $\nu_{max}$ in the algorithm as opposed to $<\nu>$.

\begin{theorem}
\label{stab:theorem}
Consider algorithm \eqref{En-Leray-POD-Weak}; define $0 \leq \epsilon \leq 1$ such that 
\begin{eqnarray}
	\max\limits_{1 \leq j \leq J} \frac{|\nu_{j} - \nu_{max}|}{ \nu_{max}} 
	= 1 - \epsilon
	\label{eqn:epsilon}
\end{eqnarray}
and assume the following condition holds for $j = 1 \ldots J$:
\begin{equation}\label{stab:assumption}
\frac{C_{b^*}^2 \, \Delta t}{\nu_{max}} \||\mathbb{S}_{R}|\|_{2}^{\frac{1}{2}}\|\nabla(\overline{u_{R}^{j,n}- <u_{R}>^{n}})\|^{2} \leq \epsilon.
\end{equation}
  Then, for any $N \geq 1$

  \begin{equation}
  \begin{aligned}
	&\frac{1}{2} \|u_{R}^{j,N}\|^{2}  
	+  \frac{\nu_{max} \Delta t}{2}\|\nabla u_{R}^{j,N}\|^{2} 
	+ \frac{\epsilon \, \nu_{max} \, \Delta t}{4} \sum_{n=0}^{N-1} \|u_{R}^{j,n+1}\|^{2}.
	\\ 
	&\leq \sum_{n=0}^{N-1} \frac{\Delta t}{\nu_{max} \epsilon}\|f_{j}^{n+1}\|^{2}_{-1} + \frac{1}{2}\|u_{R}^{0}\|^{2} +   \frac{\nu_{max} \Delta t}{2}\|\nabla u_{R}^{j,0}\|^{2} \stackrel{notation}{=} C_{stab} \, .
	\end{aligned}
	\label{eqn:theorem-stability-1}
  \end{equation}

  \begin{proof}
  	Setting $\varphi = u_{R}^{j,n+1}$ and using the skew-symmetry of the trilinear term we have
  	\begin{equation}\label{Theor1:eq1}
  	\begin{aligned}
  	&\frac{1}{2}\|u_{R}^{j,n+1}\|^{2} - \frac{1}{2}\|u_{R}^{j,n}\|^{2} + \frac{1}{2}\|u_{R}^{j,n+1} - u_{R}^{j,n}\|^{2}  + \nu_{max} \Delta t \|\nabla u_{R}^{j,n+1}\|^{2}
	\\ &+ \Delta t b^{\ast}(\overline{u_{R}^{j, n}-<u_{R}>^{n}},u_{R}^{j,n}, u_{R}^{j,n+1} - u_{R}^{j,n} )   = 
	\\&\Delta t (f_{j}^{n+1},u_{R}^{j,n+1}) - \Delta t (\nu_{j} - \nu_{max}) (\nabla u_{R}^{j,n},\nabla u_{R}^{j,n+1} ) 
  	\end{aligned}
  	\end{equation}
  	Now applying Young's inequality on the right hand side we have 
  \begin{equation}\label{Theor1:eq2}
  \begin{aligned}
  &\frac{1}{2}\|u_{R}^{j,n+1}\|^{2} - \frac{1}{2}\|u_{R}^{j,n}\|^{2} + \frac{1}{2}\|u_{R}^{j,n+1} - u_{R}^{j,n}\|^{2} + \nu_{max} \Delta t \|\nabla u_{R}^{j,n+1}\|^{2} 
  \\ 
  &+ \Delta t b^{\ast}(\overline{u_{R}^{j, n}-<u_{R}>^{n}},u_{R}^{j,n}, u_{R}^{j,n+1} - u_{R}^{j,n} )  \leq \frac{\alpha \Delta t \nu_{max} }{4}\|\nabla u_{R}^{j,n+1}\|^{2} 
  \\ 
  &+ \frac{\Delta t}{\alpha \nu_{max}} \|f_{j}^{n+1}\|^{2}_{-1} + \frac{\beta \Delta t \nu_{max} }{4}\|\nabla u_{R}^{j,n+1}\|^{2} + \frac{\Delta t(\nu_{j} - \nu_{max})^{2} }{\beta \nu_{max}} \|\nabla u_{R}^{j,n}\|^{2}.
  \end{aligned}
  \end{equation}
  Since both $\frac{\beta \Delta t \nu_{max} }{4}\|\nabla u_{R}^{j,n+1}\|^{2}$ and $\frac{\Delta t(\nu_{j} - \nu_{max})^{2} }{\beta \nu_{max}} \|\nabla u_{R}^{j,n}\|^{2} $ need to be absorbed into $\nu_{max}\Delta t \|u^{j,n+1}_{R}\|^{2}$ we minimize the quantity $\frac{\beta \Delta t \nu_{max} }{4} + \frac{\Delta t(\nu_{j} - \nu_{max})^{2} }{\beta \nu_{max}}$ by selecting $\beta = \frac{2|\nu_{j} - \nu_{max}|}{\nu_{max}} $. 
  It then follows that 
  \begin{equation}\label{Theor1:eq3}
  \begin{aligned}
  &\frac{1}{2}\|u_{R}^{j,n+1}\|^{2} - \frac{1}{2}\|u_{R}^{j,n}\|^{2} + \frac{1}{2}\|u_{R}^{j,n+1} - u_{R}^{j,n}\|^{2} + \nu_{max} \Delta t \|\nabla u_{R}^{j,n+1}\|^{2}  
  \\
  &+  \Delta tb^{\ast}(\overline{u_{R}^{j, n}-<u_{R}>^{n}},u_{R}^{j,n}, u_{R}^{j,n+1} - u_{R}^{j,n} ) \leq \frac{\alpha \Delta t \nu_{max} }{4}\|\nabla u_{R}^{j,n+1}\|^{2} 
  \\ 
  & + \frac{\Delta t}{\alpha \nu_{max}} \|f_{j}^{n+1}\|^{2}_{-1} + \frac{\Delta t|\nu_{j} - \nu_{max}| }{2}\|\nabla u_{R}^{j,n+1}\|^{2} + \frac{\Delta t|\nu_{j} - \nu_{max}| }{2} \|\nabla u_{R}^{j,n}\|^{2}.
  \end{aligned}
  \end{equation}
  Next we bound the trilinear term using \eqref{In1} and \eqref{POD:inveq}, obtaining
	\begin{equation}\label{Theor1:eq4}
	\begin{aligned}
	-\Delta t b^{\ast}&(\overline{u_{R}^{j, n}-<u_{R}
		>^{n}},u_{R}^{j,n}, u_{R}^{j,n+1} - u_{R}^{j,n} ) \\
	&\leq C_{b^\ast} \Delta t \|\overline{u_{R}^{j, n}-<u_{R}>^{n}}\| \| \nabla u_{R}^{j,n}\|\left(\|\nabla (u_{R}^{j,n+1} - u_{R}^{j,n}) \| \|u_{R}^{j,n+1} - u_{R}^{j,n}\|  \right)^{\frac{1}{2}} \\
	&\leq C_{b^\ast} \Delta t \||\mathbb{S}_{R}|\|^{\frac{1}{4}}_{2} \|\overline{u_{R}^{j, n}-<u_{R}
		>^{n}}\| \| \nabla u_{R}^{j,n}\| \|u_{R}^{j,n+1} - u_{R}^{j,n}\| \, .
	\end{aligned}
	\end{equation}
	Then using Young's inequality we obtain
	\begin{equation}\label{Theor1:eq5}
	\begin{aligned}
	-\Delta t b^{\ast}&(\overline{u_{R}^{j, n}-<u_{R}>^{n}},u_{R}^{j,n}, u_{R}^{j,n+1} - u_{R}^{j,n} ) \\
	&\leq \frac{C_{b^*}^2 \, \Delta t^{2}}{2} \||\mathbb{S}_{R}|\|^{\frac{1}{2}}_{2} \|\overline{u_{R}^{j, n}-<u_{R}
		>^{n}}\|^{2} \| \nabla u_{R}^{j,n}\|^{2} + \frac{1}{2}\|u_{R}^{j,n+1} - u_{R}^{j,n}\|^{2} \, .	 	
	\end{aligned}
	\end{equation}
	Combining like terms we then have
\begin{equation}\label{Theor1:eq6}
\begin{aligned}
&\frac{1}{2}\|u_{R}^{j,n+1}\|^{2} - \frac{1}{2}\|u_{R}^{j,n}\|^{2}  + \nu_{max} \Delta t\left(1 - \frac{\alpha}{4} - \frac{|\nu_{j} - \nu_{max}|}{2 \nu_{max}}\right)\|\nabla u_{R}^{j,n+1}\|^{2}  
\\
&\leq \frac{\Delta t}{\alpha \nu_{max}} \|f_{j}^{n+1}\|^{2}_{-1} + \frac{C_{b^*}^2 \, \Delta t^{2}}{2}  \||\mathbb{S}_{R}|\|^{\frac{1}{2}}_{2} \|\overline{u_{R}^{j, n}-<u_{R}>^{n}}\|^{2} \| \nabla u_{R}^{j,n}\|^{2} 
\\
&+ \frac{\Delta t|\nu_{j} - \nu_{max}| }{2} \|\nabla u_{R}^{j,n}\|^{2}  .
\end{aligned}
\end{equation}
Rearranging terms it follows that 
\begin{equation}\label{Theor1:eq7}
\begin{aligned}
&\frac{1}{2}\|u_{R}^{j,n+1}\|^{2} - \frac{1}{2}\|u_{R}^{j,n}\|^{2}  + \nu_{max} \Delta t \biggr( \big(1 - \frac{\alpha}{4} - \frac{|\nu_{j} - \nu_{max}|}{2 \nu_{max}}\big)\|\nabla u_{R}^{j,n+1}\|^{2}
\\ &- (\frac{|\nu_{j} - \nu_{max}|}{2 \nu_{max}} + \frac{C_{b^*}^2 \, \Delta t}{2 \, \nu_{max}} \||\mathbb{S}_{R}|\|^{\frac{1}{2}}_{2} \|\overline{u_{R}^{j, n}-<u_{R} >^{n}}\|^{2}) \| \nabla u_{R}^{j,n}\|^{2}  \biggr) \leq \frac{\Delta t}{\alpha \nu_{max}} \|f_{j}^{n+1}\|^{2}_{-1}.
\end{aligned}
\end{equation}
Using the fact that$\max\limits_{1 \leq j \leq J} \frac{|\nu_{j} - \nu_{max}|}{ \nu_{max}} = 1 - \epsilon$ for $0 \leq \epsilon \leq 1$ and taking $\alpha = \epsilon$ we have 
 \begin{equation}\label{Theor1:eq8}
\begin{aligned}
&\frac{1}{2}\|u_{R}^{j,n+1}\|^{2} - \frac{1}{2}\|u_{R}^{j,n}\|^{2}  + \nu_{max} \Delta t \biggr( \left( \frac{1}{2} + \frac{\epsilon}{4} \right) \|\nabla u_{R}^{j,n+1}\|^{2}
\\ &- (\frac{1}{2}  - \frac{\epsilon}{2} + \frac{C_{b^*}^2 \, \Delta t}{2 \, \nu_{max}} \||\mathbb{S}_{R}|\|^{\frac{1}{2}}_{2} \|\overline{u_{R}^{j, n}-<u_{R} >^{n}}\|^{2}) \| \nabla u_{R}^{j,n}\|^{2}  \biggr) \leq \frac{ \Delta t}{\nu_{max} \epsilon} \|f_{j}^{n+1}\|^{2}_{-1}.
\end{aligned}
\end{equation}
Now using assumption \eqref{stab:assumption}, \eqref{Theor1:eq8} we have
 \begin{equation}\label{Theor1:eq9}
\begin{aligned}
&\frac{1}{2}\|u_{R}^{j,n+1}\|^{2} 
- \frac{1}{2}\|u_{R}^{j,n}\|^{2}  
+ \nu_{max} \Delta t \biggr( 
	\frac{1}{2}\|\nabla u_{R}^{j,n+1}\|^{2} - \frac{1}{2}   \| \nabla u_{R}^{j,n}\|^{2}  
\biggr) 
\\
& + \nu_{max} \Delta t \frac{\epsilon}{4} \, \|u_{R}^{j,n+1}\|^{2} \leq \frac{ \Delta t}{\nu_{max} \epsilon} \|f_{j}^{n+1}\|^{2}_{-1}.
\end{aligned}
\end{equation}

Summing up~\eqref{Theor1:eq9} from $0$ to $N-1$ yields~\eqref{eqn:theorem-stability-1}.
	\end{proof}
\end{theorem}
 \\
 \begin{remark}
 The term $\epsilon$ in the above theorem measures the relative uncertainty present in the viscosities. In practice, the amount of uncertainty present in the viscosities can be one or two orders of magnitude. In this case $\epsilon \approx \mathcal{O}(10^{-1})$ or $\mathcal{O}(10^{-2})$.
 \end{remark}

\section{Error analysis}
\label{err_analysis}

We next provide an error analysis for Leray ensemble-POD solutions. 
First, we present several results obtained in \cite{GJS17}, which we use in the analysis.
We also use the following notation:

\begin{definition}[Generic Constant $C$]
Let $C$ be a generic constant that can depend on $f,u^{j}$, but not on $h, \Delta t, R, \lambda_i, \epsilon, \nu_{max}, \delta, C_{stab}, C_{b^{*}}$.
\end{definition}

%We first define the $L^2(\Omega)$ projection operator $P_R$: $L^2(\Omega) \rightarrow X_R$ by
%\begin{equation*}%\label{eq:def_proj}
%(u-P_R u , \varphi)=0\qquad \forall \varphi \in X_R.
%\end{equation*}

The following lemma is similar to Lemma 5.1 in~\cite{GJS17}.

\begin{lemma} \label{lm:L2err} 
	{\rm[$L^2(\Omega)$ norm of the error between snapshots and their projections onto the POD space]}  We have
	\begin{equation*}%\label{errL2_1}
	\frac{1}{J_S(N_S+1)} \sum_{j=1}^{J_S} \sum_{m=0}^{N_S}\Big \|  u_{h,S}^{j,m}-\sum_{i=1}^R (u_{h,S}^{j,m}, \varphi_i)\varphi_i\Big \| ^2 = \sum_{i=R+1}^{J_S(N_S+1)} \lambda_i 
	\end{equation*}
	and thus for $j=1,\ldots,J_S$,
	\begin{equation*}%\label{errL2_2}
	\frac{1}{N_S+1} \sum_{m=0}^{N_S}\Big \|  u_{h,S}^{j,m}-\sum_{i=1}^R (u_{h,S}^{j,m}, \varphi_i)\varphi_i\Big \| ^2 \leq J_S\sum_{i=R+1}^{J_S(N_S+1)} \lambda_i .
	\end{equation*}
\end{lemma}

The following lemma is similar to Lemma 5.2 in~\cite{GJS17}.

\begin{lemma}\label{lm:H1err}{\rm [$H^1(\Omega)$ norm of the error between snapshots and their projections in the POD space.]}  We have
	\begin{equation*}%\label{errH1_1}
	\frac{1}{J_S(N_S+1)} \sum_{j=1}^{J_S} \sum_{m=0}^{N_S}\Big \| \nabla \Big( u_{h,S}^{j,m}-\sum_{i=1}^R (u_{h,S}^{j,m}, \varphi_i)\varphi_i\Big)\Big\| ^2 
	=\sum_{i=R+1}^{J_S(N_S+1)} \lambda_i  \|  \nabla \varphi_i\|^2
	\end{equation*}
	and thus, for $j=1,\ldots,J_S$,
	\begin{equation*}%\label{errH1_2}
	\frac{1}{N_S+1} \sum_{m=0}^{N_S} \Big\| \nabla\Big( u_{h,S}^{j,m}-\sum_{i=1}^R (u_{h,S}^{j,m}, \varphi_i)\varphi_i \Big)\Big\| ^2 \leq J_S\sum_{i=R+1}^{J_S(N_S+1)} \lambda_i  \|  \nabla \varphi_i\|^2.
	\end{equation*}
\end{lemma}

The following lemma is similar to Lemma 5.3 in~\cite{GJS17} (see also Lemma 3.3 in~\cite{IW14}).

\begin{lemma}\label{lm:Projerr}{\rm [Error in the projection onto the POD space]}
	Consider the partition $0=t_0<t_1< \cdots < t_{N_S} = T$ used in Section \ref{POD_sec}. 
	For any $u \in H^{1}(0,T;[H^{s+1}(\Omega)]^d)$, let $u^{m}=u(\cdot, t_m)$. Then, the error in the projection onto the POD space $X_R$ satisfies the estimates
	\begin{equation*}%\label{pel2}
	\begin{aligned}
	\frac{1}{N_S+1}& \sum_{m=0}^{N_S} \|  u^{j,m}-P_R u^{j,m}\| ^2
	\leq 
	C\left( h^{2s+2} + \Delta t^4  \right)
	+ J_S\sum_{i=R+1}^{J_S(N_S+1)} \lambda_i 
	\\\\
	\frac{1}{N_S+1} &\sum_{m=0}^{N_S} \|  \nabla \left(u^{j,m}-P_R u^{j,m}\right)\| ^2
	\\&
	\leq (C+h^2 \|\hspace{-1pt}|  {\mathbb S}_R \|\hspace{-1pt}| _2 ) h^{2s} + (C+\|\hspace{-1pt}|  {\mathbb S}_R \|\hspace{-1pt}| _2 )\Delta t^4  
	+ J_S\sum_{i=R+1}^{J_S(N_S+1)} \|  \nabla \varphi_i\| ^2\lambda_i . 
	\end{aligned}
	\end{equation*}
\end{lemma}

We assume the following estimates are also valid, as done in \cite{IW14}.
\begin{assumption}\label{assumption1}
	Consider the partition $0=t_0<t_1< \cdots < t_{N_S} = T$ used in Section \ref{POD_sec}. 
	For any $u \in H^{1}(0,T;[H^{s+1}(\Omega)]^d)$, let $u^{m}=u(\cdot, t_m)$. Then, the error in the projection onto the POD space $X_R$ satisfies the estimates
	\begin{equation*}%\label{pel2}
	\begin{aligned}
	&\|  u^{j,m}-P_R u^{j,m} \| ^2 \leq C\left( h^{2s+2} + \Delta t^4  \right) + J_S\sum_{i=R+1}^{J_S(N_S+1)} \lambda_i 
	\\\\
	&\|\nabla \left(u^{j,m}-P_R u^{j,m} \right)\|^2
	\\&\leq (C+h^2 \|\hspace{-1pt}|  {\mathbb S}_R \|\hspace{-1pt}| _2 ) h^{2s} + (C+\|\hspace{-1pt}|  {\mathbb S}_R \|\hspace{-1pt}| _2 )\Delta t^4  
	+ J_S\sum_{i=R+1}^{J_S(N_S+1)} \|  \nabla \varphi_i\| ^2\lambda_i . 
	\end{aligned}
	\end{equation*}

\end{assumption}
{
Next we need to make an assumption on the regularity of $u^{j,m}_{R}$ in order to establish an estimate for the ROM filtering error. We note that this assumption is consistent with our regularity assumption \ref{assumption:reg}.
\\
\begin{assumption}\label{assumption:leray}
We assume that $\Delta u^{j,m}_{R} \in L^{2}$
\end{assumption}
}
\\\\
We now state an estimate for the ROM filtering error which is a simple extension of Lemma 4.3 in \cite{XIE201812}.

\begin{lemma} \label{lm:ROMfiltering}{\rm [ROM filtering error estimates]} 
If $\Delta u_{R}^{j,m} \in L^{2}$, then the following estimate holds:
%\blue{
%Should we replace $\|u^{j,m} - \overline{u^{j,m}}\|$ with $\|u^{j,m}_{R} - \overline{u^{j,m}_{R}}\|$?
%It seems that we're using the latter in~\eqref{eqn:theorem-convergence-1}.
%}
	\begin{equation}
	\begin{aligned}
	&
	\delta^{2}\|\nabla(u_{R}^{j,m} - \overline{u_{R}^{j,m}})\|^{2} 
	+ \|u_{R}^{j,m} - \overline{u_{R}^{j,m}}\|^{2} 
	\\& \leq 
	C\biggr( 
	C\left( h^{2s+2} + \Delta t^4  \right)
	+ J_S\sum_{i=R+1}^{J_S(N_S+1)} \lambda_i  \biggr)
	\\& 
	+ C\delta^{2}  \biggr( 
	(C+h^2 \|\hspace{-1pt}|  {\mathbb S}_R \|\hspace{-1pt}| _2 ) h^{2s} 
	+ (C+\|\hspace{-1pt}|  {\mathbb S}_R \|\hspace{-1pt}| _2 )\Delta t^4  
	+ J_S\sum_{i=R+1}^{J_S(N_S+1)} \|  \nabla \varphi_i\| ^2\lambda_i 
	\biggr) 
	\\& + C\delta^{4}\|\Delta u_{R}^{j,m} \|^{2}  . 
	\end{aligned}
	\end{equation}

\end{lemma}

Lastly we state a result for the stability of the ROM filtered variables proven in Lemma 4.4 in~\cite{XIE201812}. \\

\begin{lemma}\label{lm:ROMstability}{\rm [ROM stability estimates]} For $u \in X$, we have
	\begin{equation}
	\begin{aligned}
	&\|\overline{u}\| \leq \|u\| \\
	& \| \nabla \overline{u}\| \leq |\|S_{R}\||_{2}^{\frac{1}{2}}\|u\|. \\
	\end{aligned}
	\end{equation}
	For $u \in X_{R}$, we have
	\begin{equation}
	\| \nabla \overline{u}\| \leq \|\nabla u \|.
	\end{equation}

\end{lemma}

Let $e^{j,n}=u^{j,n}-u_{R}^{j,n}$ denote the error between the true solution
and the POD approximation; then, we have the following error estimates.

\begin{theorem}
	Consider the Leray ensemble-POD algorithm and the partition $0 = t_{0} < t_{1} < \cdots < t_{N_{S}}$ used in Section \ref{POD_sec}. Suppose for any $0 \leq n \leq N_{S}$, the stability conditions from Theorem \ref{stab:theorem} and all previously stated regularity assumptions hold. Then for any $1 \leq N \leq N_{S}$, there is a positive constant $C$, such that the following bound holds: 
	
	\begin{equation}
	\begin{aligned}
	&\frac{1}{2}\|e^{j,N}\|^{2} + \frac{\nu_{max}}{2}\| \nabla e^{j,N}\|^{2} + \frac{\epsilon}{4}\nu_{max} \| \nabla e^{j,N}\|^{2}  + C\nu_{max}\Delta t \sum_{n=0}^{N-1} \| e^{j,n+1} \|^{2} 	
	\\ &\leq \exp\left(\frac{C_{b^{*}}^{4} \, C^{4} \,T}
	{\epsilon^{3}\,\nu_{max}^{3}}\right) 
	\biggr[\biggr(\frac{ C \, \nu_{max} \Delta t}{\epsilon}
	+ \frac{C \Delta t}{\epsilon} \frac{|\nu_j - \nu_{max}|^{2}}{\nu_{max}} 
	+ 2C \||\mathbb{S}_{R}|\|_{2}^{-\frac{1}{2}}
	+ \frac{C \, C_{b^*}^2 \, \Delta t}{\epsilon \, \nu_{max}}
	\\ &+ \frac{C \, C_{b^*}^2 \, C_{stab}}{2 \, \epsilon \, \nu_{max}} 
	+ \frac{2 \,C \, C_{b^*}^2 \, C^{2}_{stab}}{\epsilon^{2} \, \nu^{2}_{max}}
	+ \frac{C \, C_{b^*}^2 \, \delta}{ \epsilon \, \nu_{max}}  \biggr) \times 
	\\& \biggr((C+h^2 \|\hspace{-1pt}|  {\mathbb S}_R \|\hspace{-1pt}| _2 ) h^{2s} + (C+\|\hspace{-1pt}|  {\mathbb S}_R \|\hspace{-1pt}| _2 )\Delta t^4  
	+ J_S\sum_{i=R+1}^{J_S(N_S+1)} \|  \nabla \varphi_i\| ^2\lambda_i \biggr) 
	\\& + \frac{C\,C_{b^*}^2}{\delta \, \epsilon \, \nu_{max}} \biggr(C\left( h^{2s+2} + \Delta t^4  \right)+ J_S\sum_{i=R+1}^{J_S(N_S+1)} \lambda_i  \biggr)
	\\& + \frac{C \, \Delta t^{2}}{\epsilon}\frac{|\nu_j - \nu_{max}|^{2}}{\nu_{max}} 
	+ C \Delta t \||\mathbb{S}_{R}|\|_{2}^{-\frac{1}{2}}
	+ \frac{C h^{2s}} {d \, \epsilon \, \nu_{max}}\| |p^{j} | \|^{2}_{2,s} 
	+ \frac{C \Delta t^{2}}{\epsilon \, \nu_{max}} 
	\\& 
	+ \frac{C\,C_{b^*}^2 \, \Delta t^{2}}{\epsilon \, \nu_{max}} 
	+ \frac{C\,C_{b^*}^2 \, \Delta t \, \delta^3}{\epsilon \, \nu_{max} \, }\biggr] 
	\\&
	+ \left(1 + CN\nu_{max}\Delta t \right) \times \left(C\left( h^{2s+2} + \Delta t^4  \right) + J_S\sum_{i=R+1}^{J_S(N_S+1)} \lambda_i \right)
	\\& 
	+(\nu_{max} + \frac{\epsilon}{2}\nu_{max}) \times 
	\\&
	\left((C+h^2 \|\hspace{-1pt}|  {\mathbb S}_R \|\hspace{-1pt}| _2 ) h^{2s} + (C+\|\hspace{-1pt}|  {\mathbb S}_R \|\hspace{-1pt}| _2 )\Delta t^4  
	+ J_S\sum_{i=R+1}^{J_S(N_S+1)} \|  \nabla \varphi_i\| ^2\lambda_i \right).	
	\end{aligned}
	\end{equation}

	\begin{proof}
	
	The weak solution of the NSE $u^{j}$ satisfies
	\begin{equation}
	\label{NSE:weak}
	\begin{aligned}
	\left(\frac{u^{j,n+1} - u^{j,n}}{\Delta t},\varphi\right)& + b^{*}(u^{j,n+1},u^{j,n+1},\varphi) + \nu_{j} (\nabla u^{j,n+1}, \nabla \varphi) - (p^{j,n+1},\nabla \cdot \varphi) \\
	&= (f^{j,n+1},\varphi) + Intp(u^{j,n+1};\varphi)
	\end{aligned}
	\end{equation}
	where 
	\begin{equation}
	Intp(u^{j,n+1};\varphi) = (\frac{u^{j,n+1} - u^{j,n}}{\Delta t} - u^{j}_{t}(t^{n+1}),\varphi).
	\end{equation}
	We split the error
	\begin{equation}
	e^{j,n} = u^{j,n} - u_{R}^{j,n} = (u^{j,n} - P_{R} u^{j,n}) + (P_{R}u^{j,n} - u_{R}^{j,n}) = \eta^{j,n} + \xi_{R}^{j,n}, \qquad  j = 1,\hdots, J.
	\end{equation}
	Subtracting \eqref{En-Leray-POD-Weak} from \eqref{NSE:weak} 
	as well as adding and subtracting the terms \newline $\nu_{max}(\nabla u_{j}^{n+1},\nabla \varphi)$ and $\nu_{j} - \nu_{max} ( \nabla u_{j}^{n+1},\nabla \varphi)$ we have
	\begin{equation}
	\begin{aligned}
	(\frac{\xi_{R}^{j,n+1} - \xi_{R}^{j,n}}{\Delta t} &, \varphi) + \nu_{max}(\nabla \xi_{R}^{j,n+1},\nabla \varphi) + (\nu_{j} - \nu_{max})(\nabla(u^{j,n+1}-u^{j,n}),\nabla  \varphi)\\
	&+ (\nu_j - \nu_{max})(\nabla \xi_{R}^{j,n},\nabla \varphi) + b^{*}(u^{j,n+1},u^{j,n+1},\varphi)  \\ 
	&- b^{*}(\overline{<u_{R}
		>^{n}},u_{R}^{j,n+1},\varphi)
	- b^{*}(\overline{u_{R}^{j, n}-<u_{R}
		>^{n}},u_{R}^{j,n},\varphi) \\ 
	&- (p^{j,n+1},\nabla \cdot \varphi) \\
	& = -(\frac{\eta^{j,n+1} - \eta^{j,n}}{\Delta t}, \varphi) - \nu_{max}(\nabla \eta^{j,n+1},\nabla \varphi)  \\ &- (\nu_j - \nu_{max})(\nabla \eta^{j,n},\nabla \varphi) + Intp(u^{j,n+1};\varphi).
	\end{aligned}
	\end{equation}
	Setting $\varphi = \xi_{R}^{j,n+1}$ rearranging the nonlinear terms by adding and subtracting \newline $b^{*}(\overline{u_{R}^{j, n}-<u_{R}
		>^{n}},u_{R}^{j,n+1} ,\xi_{R}^{j,n+1})$ , and using the fact that $(\eta^{j,n+1} - \eta^{j,n}, \xi_{R}^{j,n+1}) = 0$ by the definition of the $L^{2}$ projection we have
	\begin{equation}
	\label{eq:err1}
	\begin{aligned}
	\frac{1}{\Delta t}&\left(\frac{1}{2}\|\xi_{R}^{j,n+1}\|^{2} - \frac{1}{2}\|\xi_{R}^{j,n}\|^{2} + \frac{1}{2}\|\xi_{R}^{j,n+1} - \xi_{R}^{j,n} \|^{2} \right) + \nu_{max} \|\nabla \xi_{R}^{j,n+1}\|^{2} \\
	&= - (\nu_{j} - \nu_{max})(\nabla(u^{j,n+1}-u^{j,n}),\nabla  \xi_{R}^{j,n+1})  - (\nu_j - \nu_{max})(\nabla \xi_{R}^{j,n},\nabla \xi_{R}^{j,n+1}) \\
	&  - \nu_{max}(\nabla \eta^{j,n+1},\nabla \xi_{R}^{j,n+1}) - (\nu_j - \nu_{max})(\nabla \eta^{j,n},\nabla \xi_{R}^{j,n+1})\\
	& + b^{*}(\overline{u_{R}
		^{j,n}},u_{R}^{j,n+1},\xi_{R}^{j,n+1})
	- b^{*}(\overline{u_{R}^{j, n}-<u_{R}
		>^{n}},u_{R}^{j,n+1} - u_{R}^{j,n},\xi_{R}^{j,n+1})
	\\
	&- b^{*}(u^{j,n+1},u^{j,n+1},\xi_{R}^{j,n+1}) + (p^{j,n+1},\nabla \cdot \xi_{R}^{j,n+1}) + Intp(u^{j,n+1};\xi_{R}^{j,n+1}).
	\end{aligned}
	\end{equation}
	We bound the viscous terms in a similar manner to Theorem 3.1 of \cite{GJW17}
	\begin{equation}
	\label{err1eq}
	\begin{aligned}
	-(\nu_{j} - \nu_{max})&(\nabla(u^{j,n+1} - u^{j,n}),\nabla \xi^{j,n+1}_{R}) \leq \\ &\frac{\Delta t}{4\tepsilon}\frac{|\nu_{j} - \nu_{max}|^{2}}{\nu_{max}}\left(\int_{t^{n}}^{t^{n}+1}\|\nabla u_{j,t}\|^{2}dt\right) + \tepsilon\nu_{max}\|\nabla \xi_{R}^{j,n+1}\|^{2},
	\end{aligned}
	\end{equation}
	\begin{equation}
	-\nu_{max}(\nabla \eta^{j,n+1},\nabla \xi_{R}^{j,n+1}) \leq \frac{ \nu_{max}}{4\tepsilon}\|\nabla \eta^{j,n+1}\|^{2} + \tepsilon\nu_{max}\|\nabla \xi_{R}^{j,n+1}\|^{2},
	\end{equation}	 
	\begin{equation}
	\begin{aligned}
	-(\nu_{j} - \nu_{max})&(\nabla \eta^{j,n},\nabla \xi_{R}^{j,n+1}) \leq \\ 
	&\frac{1}{4\tepsilon}\frac{|\nu_{j} - \nu_{max}|^{2}}{\nu_{max}}\|\nabla \eta^{j,n}\|^{2} + \tepsilon\nu_{max}\|\nabla \xi_{R}^{j,n+1}\|^{2},
	\end{aligned}
	\end{equation}
	\begin{equation}
	-(\nu_{j} - \nu_{max}) (\nabla \xi_{R}^{j,n},\nabla \xi_{R}^{j,n+1}) 
	\leq \frac{|\nu_{j} - \nu_{max}|}{2}\| \nabla \xi_{R}^{j,n}\|^{2} 
	+ \frac{|\nu_{j} - \nu_{max}|}{2}\| \nabla \xi_{R}^{j,n+1}\|^{2}.
	\end{equation}
	We next rewrite the second nonlinear term on the right hand side of \eqref{eq:err1}.
	\begin{equation}
	\begin{aligned}
	b^{*}&(\overline{u_{R}^{j, n}-<u_R>^n},u_{R}^{j,n+1}-u_{R}^{j,n},\xi_{R}^{j,n+1})\\
	&=-b^{*}(\overline{u_{R}^{j, n}-<u_{R}
		>^{n}},e^{j,n+1}-e^{j,n},\xi_{R}^{j,n+1})\\
	&\quad+b^{*}(\overline{u_{R}^{j, n}-<u_{R}
		>^{n}}
	,u^{j,n+1}-u^{j,n},\xi_{R}^{j,n+1})\\
	&=-b^{*}(\overline{u_{R}^{j, n}-<u_{R}
		>^{n}},\eta^{j,n+1},\xi_{R}^{j,n+1})\\
	&\quad+b^{*}(\overline{u_{R}^{j, n}-<u_{R}
		>^{n}},\eta
	^{j,n},\xi_{R}^{j,n+1})\\
	&\quad+b^{*}(\overline{u_{R}^{j, n}-<u_{R}
		>^{n}},\xi_{R}^{j,n},\xi_{R}^{j,n+1})\\
	&\quad+b^{*}(\overline{u_{R}^{j, n}-<u_{R}
		>^{n}},u
	^{j,n+1}-u^{j,n},\xi_{R}^{j,n+1})\text{ .}
	\end{aligned}
	\end{equation}
	As done in Theorem 3.1 of \cite{GJW17} using Young's inequality, \eqref{In1}, \eqref{In2}, and \eqref{POD:inveq} we derive the estimates
	
	\begin{equation}
	\begin{aligned}
	-b^{*}(&\overline{u_{R}^{j, n}-<u_{R}
		>^{n}},\eta^{j,n+1},\xi_{R}^{j,n+1}) \leq \\ &\frac{C_{b^*}^{2}\nu_{max}^{-1}}{4\tepsilon}\|\nabla(\overline{u_{R}^{j, n}-<u_{R}
		>^{n}})\|^{2}\|\nabla \eta^{j,n+1}\|^{2} + \tepsilon\nu_{max}\|\nabla \xi_{R}^{j,n+1}\|^{2},
	\end{aligned}
	\end{equation}
	
		\begin{equation}
		\begin{aligned}
		b^{*}(&\overline{u_{R}^{j, n}-<u_{R}
			>^{n}},\eta^{j,n},\xi_{R}^{j,n+1}) \leq \\ &\frac{C_{b^*}^{2}\nu_{max}^{-1}}{4\tepsilon}\|\nabla(\overline{u_{R}^{j, n}-<u_{R}
			>^{n}})\|^{2}\|\nabla \eta^{j,n}\|^{2} + \tepsilon\nu_{max}\|\nabla \xi_{R}^{j,n+1}\|^{2},
		\end{aligned}
		\end{equation}
	 
	 \begin{equation}
	 \begin{aligned}
	 &b^{*}(\overline{u_{R}^{j, n}-<u_{R}
	 	>^{n}},u^{j,n+1} - u^{j,n},\xi_{R}^{j,n+1}) \leq \\ &\frac{C \,C_{b^*}^{2}\nu_{max}^{-1}}{4\tepsilon} \Delta t \|\nabla(\overline{u_{R}^{j, n}-<u_{R}
	 	>^{n}})\|^{2} + \tepsilon\nu_{max}\|\nabla \xi_{R}^{j,n+1}\|^{2}.
	 \end{aligned}
	 \end{equation}
	\noindent  By skew-symmetry, inequality \eqref{In2} and the inverse inequality \eqref{POD:inveq}, we have
	\begin{equation}
	\begin{aligned}
	&b^{\ast}(\overline{u_{R}^{j, n}-<u_{R}
		>^{n}},\xi_{R}^{j,n},\xi_{R}^{j,n+1})\\
	&\leq C_{b^*} \, \| \nabla (\overline{u_{R}^{j, n}-<u_{R}
		>^{n}})\| \| \nabla\xi_{R}^{j,n}\| \sqrt{\| \xi_{R}^{j,n+1}-\xi_{R}
		^{j,n}\| \|  \nabla (\xi_{R}^{j,n+1}-\xi_{R}
		^{j,n})\| }\\
		&\leq C_{b^*} \, \|\hspace{-1pt}|  {\mathbb S}_R \|\hspace{-1pt}| _2^{1/4}| \nabla (\overline{u_{R}^{j, n}-<u_{R}
			>^{n}})\| \| \nabla\xi_{R}^{j,n}\|\| \xi_{R}^{j,n+1}-\xi_{R}
		^{j,n}\| \\
	&
		\leq \frac{1}{2\Delta t}\| \xi_{R}^{j,n+1}-\xi_{R}^{j,n}\| ^{2}+\left(
	\frac{C_{b^*}^2 \Delta t}{2} \|\hspace{-1pt}|  {\mathbb S}_R \|\hspace{-1pt}| _2^{\frac{1}{2}}\| \nabla (\overline{u_{R}^{j, n}-<u_{R}
		>^{n}})\| ^{2}\right)  \| \nabla
	\xi_{R}^{j,n}\| ^{2}.
	\end{aligned}
	\label{eqn:5.18}
	\end{equation}
	Bounding the other two nonlinear terms we add and subtract the terms \newline $b^{*}(u^{j,n},u^{j,n+1},\xi_{R}^{j,n+1})$ and $b^{*}(\overline{u^{j,n}_R},u^{j,n+1},\xi_{R}^{j,n+1})$. It then follows from \eqref{In1}    
	\begin{equation}
	\begin{aligned}
	&- b^{*}(u^{j,n+1},u^{j,n+1},\xi_{R}^{j,n+1}) +b^{*}(\overline{u_{R}
		^{j,n}},u_{R}^{j,n+1},\xi_{R}^{j,n+1}) \\
	& = -b^{*}(u^{j,n} - \overline{u_{R}
		^{j,n}} , u^{j,n+1}, \xi_{R}^{j,n+1}) -b^{*}(\overline{u_{R}
		^{j,n}},\eta^{j,n+1},\xi_{R}^{j,n+1}) \\ &-b^{*}(u^{j,n+1} - u^{j,n}, u^{j,n+1}, \xi_{R}^{j,n+1}).
	\end{aligned}
	\end{equation}
	
	Now by Young's inequality, \eqref{In2}, the stability analysis, i.e. 
	
	$\|\overline{u^{j,n}_{R}}\|^{2} \leq C_{stab}$, and the assumption $u^{j} \in L^{\infty}(0,T,H^{1}(\Omega))$  we have
	
	\begin{equation}
	\begin{aligned}
	b^{*}(\overline{u_{R}^{j,n}},\eta^{j,n+1},\xi_{R}^{j,n+1}) &\leq C_{b^*}\|\nabla \overline{u^{j,n}_{R}}\|^{\frac{1}{2}}\| \overline{u^{j,n}_{R}}\|^{\frac{1}{2}}\|\nabla \eta^{j,n+1}\|\|\nabla \xi_{R}^{j,n+1}\| \\ 
	&\leq \frac{C_{stab}C_{b^*}^{2}}{4\tepsilon}\nu_{max}^{-1}  \|\nabla \overline{u_{R}^{j,n}}\| \|\nabla \eta^{j,n+1}\|^{2} + 
	\tepsilon\nu_{max}\|\nabla \xi_{R}^{j,n+1}\|^{2},
	\end{aligned}
	\end{equation}
	as well as
	\begin{equation}
	\begin{aligned}
	b^{*}(u^{j,n+1} - u^{j,n}, u^{j,n+1}, &\xi_{R}^{j,n+1}) \leq \frac{C \, C_{b^*}^{2}\Delta t}{4\tepsilon}\nu_{max}^{-1} + \tepsilon\nu_{max}\|\nabla \xi_{R}^{j,n+1}\|^{2}.
	\end{aligned}
	\end{equation}
	We can then rewrite the term 
	\begin{equation}
	\begin{aligned}
	-b^{*}(u^{j,n} - \overline{u_{R}
		^{j,n}} , &u^{j,n+1}, \xi_{R}^{j,n+1}) =\\
		 &-b^{*}(e^{j,n} , u^{j,n+1}, \xi_{R}^{j,n+1}) -  b^{*}(u_{R}^{j,n} - \overline{u_{R}^{j,n}}, u^{j,n+1}, \xi_{R}^{j,n+1}). 
	\end{aligned}
	\end{equation}
	Bounding the second term
	\begin{equation}
	\begin{aligned}
	-b^{*}(u_{R}^{j,n} - &\overline{u_{R}
		^{j,n}} , u^{j,n+1}, \xi_{R}^{j,n+1}) \\ &\leq C_{b^*}\|u_{R}^{j,n} - \overline{u_{R}^{j,n}}\|^{\frac{1}{2}} \|\nabla (u_{R}^{j,n} - \overline{u_{R}^{j,n}})\|^{\frac{1}{2}}\|\nabla u^{j,n+1}\|\| \nabla \xi^{j,n+1}_{R}\|
	\\ & \leq \frac{C_{b^*}^{2}}{4\tepsilon}\nu_{max}^{-1}\|u_{R}^{j,n} - \overline{u_{R}^{j,n}}\| \|\nabla (u_{R}^{j,n} - \overline{u_{R}^{j,n}})\|\|\nabla u^{j,n+1}\|^{2} + \tepsilon\nu_{max}\|\nabla \xi^{j,n+1}_{R}\|^{2}.
	\end{aligned}
	\end{equation}
	Then after decomposing $e^{j,n} = \eta^{j,n} + \xi_{R}^{j,n}$ again using Young's inequality
 and the assumption $u^{j} \in L^{\infty}(0,T,H^{1}(\Omega))$
	\begin{equation}
	-b^{*}(\eta^{j,n} , u^{j,n+1}, \xi_{R}^{j,n+1}) 
	\leq 
	\frac{C C_{b^*}^2}{4 \, \tepsilon} \nu_{max}^{-1} \, \| \nabla \eta^{j,n}\|^{2} 
	+ \tepsilon\nu_{max}\|\nabla \xi^{j,n+1}_{R}\|^{2}
	\end{equation}
	and
	\begin{equation}
	\begin{aligned}
	-b^{*}(\xi_{R}^{j,n} , u^{j,n+1}, \xi_{R}^{j,n+1}) 
	&\leq C_{b^*} \|\nabla \xi_{R}^{j,n}\|^{\frac{1}{2}}\|\xi_{R}^{j,n}\|^{\frac{1}{2}}\|\nabla u^{j,n+1}\|\|\nabla \xi_{R}^{j,n+1}\| 
	\\ &\leq C_{b^*} C 
	\left(
		\alpha\|\nabla \xi_{R}^{j,n+1}\|^{2} 
		+ \frac{1}{4 \alpha} \|\nabla \xi_{R}^{j,n}\|\|\xi_{R}^{j,n}\|
	\right). 
	\\ &\leq C_{b^*} C 
	\left(
		\alpha\|\nabla \xi_{R}^{j,n+1}\|^{2} 
		+ \frac{1}{4 \alpha} 
		\left(
			\beta \| \nabla \xi_{R}^{j,n} \|^{2} 
			+ \frac{1}{\beta}\| \xi_{R}^{j,n} \|^{2}
		\right)
	\right)
 	\\ & {=} \, \tepsilon\nu_{max} \|\nabla \xi_{R}^{j,n+1}\|^{2} 
	+ \frac{{13} \, \tepsilon}{4}\nu_{max}\| \nabla \xi_{R}^{j,n} \|^{2} 
	+ \frac{{ C_{b^*}^4 \, C^4} }{52 \, \nu_{max}^{3}\tepsilon^{3}}\| \xi_{R}^{j,n} \|^{2}.
	\end{aligned}
	\end{equation}

	For the pressure term since $\xi_{j,r}^{n+1} \in X^{R} \subset V^{h}$ it follows for $q_{h} \in Q^{h}$
	\begin{equation}
	\begin{aligned}
	(p_{j}^{n+1}, \nabla \cdot \xi^{j,n+1}_{R}) =& (p^{j,n+1}- q_{h}^{n+1}, \nabla \cdot \xi_{R}^{j,n+1})\\
	&\leq \tepsilon\nu_{max}\|\nabla \xi_{R}^{n+1}\|^{2} + \frac{\nu_{max}^{-1}}{4 d\tepsilon}\|p^{j,n+1} - q_{h}^{j,n+1} \|^{2}.
	\end{aligned}
	\end{equation}
	For the last term we have 
	\begin{equation}
	\label{errlast}
	\begin{aligned}
	Intp(u^{j,n+1},\xi_{R}^{j,n+1}) 
	&\leq \|\frac{u^{j,n+1}-u^{j,n}}{\Delta t} - u^{j}_{t}(t^{n+1})\|\|\nabla \xi_{R}^{j,n+1}\| \\
	& \leq \frac{C \Delta t}{4\tepsilon} \nu_{max}^{-1} t
	+ \tepsilon\nu_{max}\|\nabla \xi_{R}^{j,n+1}\|^{2}.
	\end{aligned}
	\end{equation}
	
	Now combining \eqref{err1eq} - \eqref{errlast}, \eqref{eq:err1} becomes
	\begin{equation}\label{comb1}
	\begin{aligned}
	& \frac{1}{\Delta t}\left(\frac{1}{2}\|\xi_{R}^{j,n+1}\|^{2} - \frac{1}{2}\|\xi_{R}^{j,n}\|^{2}  \right)
	\\
	& + \nu_{max}\|\nabla \xi_{R}^{j,n+1}\|^{2} - \frac{C_{b^*}^2 \Delta t}{2} \|\hspace{-1pt}|  {\mathbb S}_R \|| _2^{\frac{1}{2}}\| \nabla (\overline{u_{R}^{j, n}-<u_{R}>^{n}})\| ^{2}\    \| \nabla \xi_{R}^{j,n}\| ^{2} 
	\\
	&
	  - \frac{13 \tilde{\epsilon}}{4}\nu_{max}\|\nabla \xi_{R}^{j,n+1}\|^{2} 
	  - \frac{13 \tilde{\epsilon}}{4} \, \nu_{max}\|\nabla \xi_{R}^{j,n}\|^{2} 
	\\
	&
		- \frac{|\nu_{j} - \nu_{max}|}{2}\|\nabla \xi_{R}^{j,n+1}\|^{2}  
		- \frac{|\nu_{j} - \nu_{max}|}{2}\|\nabla \xi_{R}^{j,n}\|^{2} 
	\\ 
	&
	\leq 
	{
	\frac{C_{b^*}^4 \, C^4}{52 \, \nu_{max}^{3} \, \tepsilon^{3}}\| \xi_{R}^{j,n} \|^{2} 
	}
	+ \frac{C\Delta t}{4 \, \tepsilon} \frac{|\nu_{j} - \nu_{max}|^{2}}{\nu_{max}}  
	+ \frac{\nu_{max}}{4 \, \tepsilon} \, \|\nabla \eta^{j,n+1}\|^{2}   
	\\ 
	&
	+ \frac{C C_{b^*}^2 \, \Delta t}{4 \, \tepsilon \, \nu_{max}} \, \|\nabla(\overline{u_{R}^{j, n}-<u_{R}>^{n}})\|^{2}
	\\ 
	&
	+ \frac{1}{4 \, d \, \tepsilon \, \nu_{max}} \, \|p^{j,n+1} - q_{h}^{n+1} \|^{2} 
	+ \frac{C \Delta t}{4 \, \tepsilon \, \nu_{max}} \,   
	\\ 
	&
	+  \frac{C C_{b^*}^2 \, \Delta t}{4 \, \tepsilon \nu_{max}} \, 
	+ \frac{C_{b^*}^2}{4 \, \tepsilon \nu_{max}} \, \|u_{R}^{j,n} - \overline{u_{R}^{j,n}}\| \|\nabla (u_{R}^{j,n} - \overline{u_{R}^{j,n}})\|\|\nabla u^{j,n+1}\|^{2}
	\\
	&
	+ \frac{1}{4 \, \tepsilon} \, \frac{|\nu_{j} - \nu_{max}|^{2}}{\nu_{max}} \, \|\nabla \eta^{j,n}\|^{2}
	+ \frac{C_{b^*}^{2}}{4 \, \tepsilon \, \nu_{max}} \, \|\nabla(\overline{u_{R}^{j, n}-<u_{R}>^{n}})\|^{2} \, \|\nabla \eta^{j,n}\|^{2}
	\\
	&
	+ \frac{C_{b^*}^{2}}{4 \, \tepsilon \, \nu_{max}} \, \|\nabla(\overline{u_{R}^{j, n}-<u_{R}>^{n}})\|^{2} \, \|\nabla \eta^{j,n+1}\|^{2}
	+ \frac{C_{b^*}^{2} \, {C_{stab}} }{4 \, \tepsilon \, \nu_{max}} \, \|\nabla \overline{u_{R}^{j,n}}\| \|\nabla \eta^{j,n+1}\|^{2} 
	\\ 
	&
	+ \frac{C_{b^*}^{2} \, C }{4 \, \tepsilon \, \nu_{max}} \, \|\nabla \eta^{j,n}\|^{2} 
	\, .
	\end{aligned}
	\end{equation}

The terms on the LHS of~\eqref{comb1} (except first) can be rearranged as follows:
\begin{equation}\label{eqn:traian-1}
\begin{aligned}
	\left(
		\nu_{max}
		- \frac{13 \tilde{\epsilon}}{4} \nu_{max}
		- \frac{|\nu_{j} - \nu_{max}|}{2}
	\right)		
	&\|\nabla \xi_{R}^{j,n+1}\|^{2}
	 \\
	-
	\left(
		\frac{13 \tepsilon}{4} \, \nu_{max}
		+ \frac{|\nu_{j} - \nu_{max}|}{2}
		+ \frac{C_{b^*}^2 \Delta t}{2} \|\hspace{-1pt}|  {\mathbb S}_R \|| _2^{\frac{1}{2}}\| \nabla (\overline{u_{R}^{j, n}-<u_{R}>^{n}})\| ^{2}
	\right)		
	&\|\nabla \xi_{R}^{j,n}\|^{2}	\, .
	\end{aligned}
\end{equation}
%\begin{eqnarray}
%	\left(
%		\nu_{max}
%		- \frac{13 \tilde{\epsilon}}{4} \nu_{max}
%		- \frac{|\nu_{j} - \nu_{max}|}{2}
%	\right)		
%	\|\nabla \xi_{R}^{j,n+1}\|^{2}
%	\nonumber \\
%	-
%	\left(
%		\frac{13 \tepsilon}{4} \, \nu_{max}
%		+ \frac{|\nu_{j} - \nu_{max}|}{2}
%		+ \frac{C_{b^*}^2 \Delta t}{2} \|\hspace{-1pt}|  {\mathbb S}_R \|| _2^{\frac{1}{2}}\| \nabla (\overline{u_{R}^{j, n}-<u_{R}>^{n}})\| ^{2}
%	\right)		
%	\|\nabla \xi_{R}^{j,n}\|^{2}	\, .
%	\label{eqn:traian-1}	
%\end{eqnarray}
Choosing $\tilde{\epsilon} = \frac{\epsilon}{13}$ and using~\eqref{eqn:epsilon} in~\eqref{eqn:traian-1}, \eqref{comb1} yields

	\begin{equation}
	\label{comb2}
	\begin{aligned}
	&\frac{1}{\Delta t}\left(\frac{1}{2}\|\xi_{R}^{j,n+1}\|^{2} - \frac{1}{2}\|\xi_{R}^{j,n}\|^{2} \right) 
	\\ &
	+ \left(
		\frac{\nu_{max}}{2} + \frac{\epsilon \, \nu_{max}}{4}
	\right) 
	\left(
		\| \nabla \xi_{R}^{j,n+1} \|^{2} - \| \nabla \xi_{R}^{j,n} \|^{2}
	\right)
	\\ & + 
	{\nu_{max}}\biggl(\frac{\epsilon}{2} - \frac{C_{b^*}^2 \, \Delta t}{2 \nu_{max}} \||\mathbb{S}_{R}|\|^{\frac{1}{2}}_{2} \|\nabla(\overline{u_{R}^{j, n}-<u_{R} >^{n}})\|^{2} \biggr)\| 
	\nabla \xi_{R}^{j,n} \|^{2}
	\\ &
	\leq 
	{
	\frac{13^3 \, C_{b^*}^4 \, C^4}{52 \, \nu_{max}^{3} \, \epsilon^{3}}\| \xi_{R}^{j,n} \|^{2} 
	}
	+ \frac{13 \, C\, \Delta t}{4 \, \epsilon} \frac{|\nu_{j} - \nu_{max}|^{2}}{\nu_{max}} \,  
	+ \frac{13 \, \nu_{max}}{4 \, \epsilon} \, \|\nabla \eta^{j,n+1}\|^{2}   
	\\ 
	&
	+ \frac{13 \, C \, C_{b^*}^2 \, \Delta t}{4 \, \epsilon \, \nu_{max}} \, \|\nabla(\overline{u_{R}^{j, n}-<u_{R}>^{n}})\|^{2} \, 
	\\ 
	&
	+ \frac{13}{4 \, d \, \epsilon \, \nu_{max}} \, \|p^{j,n+1} - q_{h}^{n+1} \|^{2} 
	+ \frac{13\, C \, \Delta t}{4 \, \epsilon \, \nu_{max}} \, 
	\\ 
	&
	+  \frac{13 \, C \, C_{b^*}^2 \, \Delta t}{4 \, \epsilon \nu_{max}} \,  
	+ \frac{13 \, C_{b^*}^2}{4 \, \epsilon \nu_{max}} \, \|u_{R}^{j,n} - \overline{u_{R}^{j,n}}\| \|\nabla (u_{R}^{j,n} - \overline{u_{R}^{j,n}})\|\|\nabla u^{j,n+1}\|^{2}
	\\
	&
	+ \frac{13}{4 \, \epsilon} \, \frac{|\nu_{j} - \nu_{max}|^{2}}{\nu_{max}} \, \|\nabla \eta^{j,n}\|^{2}
	+ \frac{13 \, C_{b^*}^{2}}{4 \, \epsilon \, \nu_{max}} \, \|\nabla(\overline{u_{R}^{j, n}-<u_{R}>^{n}})\|^{2} \, \|\nabla \eta^{j,n}\|^{2}
	\\
	&
	+ \frac{13 \, C_{b^*}^{2}}{4 \, \epsilon \, \nu_{max}} \, \|\nabla(\overline{u_{R}^{j, n}-<u_{R}>^{n}})\|^{2} \, \|\nabla \eta^{j,n+1}\|^{2}
	+ \frac{13 \, C_{b^*}^{2} \, {C_{stab}} }{4 \, \epsilon \, \nu_{max}} \, \|\nabla \overline{u_{R}^{j,n}}\| \|\nabla \eta^{j,n+1}\|^{2} 
	\\
	&
	+ \frac{13 \, C_{b^*}^{2} \, {C}}{4 \, \epsilon \, \nu_{max}} \, \|\nabla \eta^{j,n}\|^{2} 
	\, .
	\end{aligned}
	\end{equation}

\clearpage

	It follows from the stability condition \eqref{stab:assumption} that 
	\begin{equation}
	\label{stab:cond}
	{\nu_{max}} \left(\frac{\epsilon}{2} - \frac{{C_{b^*}^2} \, \Delta t}{2 \nu_{max}} \||\mathbb{S}_{R}|\|^{\frac{1}{2}}_{2} \|\nabla(\overline{u_{R}^{j, n}-<u_{R} >^{n}})\|^{2}\right) \geq {C \nu_{max}} \geq 0
	\end{equation}
		
	Now we use \eqref{stab:cond}, sum \eqref{comb2} from $n=0$ to $N-1$, multiply both sides by $\Delta t$, and absorb constants. Since $U_{R}^{j,0} = \sum_{i=1}^{R}(u^{j,0},\varphi_i)\varphi_i$, we have $\|\xi_{R}^{j,0}\|^{2} = 0$ and $\|\nabla \xi_{R}^{j,0}\|^{2} = 0$. It then follows from \eqref{comb2} that we have
	\begin{equation}\label{comb3}
	\begin{aligned}
	&\frac{1}{2}\|\xi_{R}^{j,N}\|^{2} + \frac{\nu_{max}}{2}\| \nabla \xi_{R}^{j,N}\|^{2}  + {\frac{\epsilon}{4}\nu_{max}\| \nabla \xi_{R}^{j,N}\|^{2}} +  C\nu_{max}\Delta t \sum_{n=0}^{N-1} \| \nabla \xi_{R}^{j,n+1} \|^{2} \\
	&
	\leq \Delta t \sum_{n=0}^{N-1 }\bigg \{ {\frac{C_{b^*}^4 \, C^4}{\epsilon^{3} \, \nu_{max}^{3}}\| \xi_{R}^{j,n} \|^{2}}  + \frac{C\Delta t}{{\epsilon}}\frac{|\nu_{j} - \nu_{max}|^{2}}{\nu_{max}} + \frac{C \nu_{max}}{{\epsilon}}\|\nabla \eta^{j,n+1}\|^{2}   
	\\ &+ \frac{C}{{\epsilon}} \frac{|\nu_{j} - \nu_{max}|^{2}}{\nu_{max}}\|\nabla \eta^{j,n}\|^{2} + {\frac{C\,C_{b^*}^2}{\epsilon \, \nu_{max}}} \|\nabla(\overline{u_{R}^{j, n}-<u_{R} >^{n}})\|^{2}\|\nabla \eta^{j,n+1}\|^{2}
	\\ &+ {\frac{C\,C_{b^*}^2}{\epsilon \, \nu_{max}}}\|\nabla(\overline{u_{R}^{j, n}-<u_{R}>^{n}})\|^{2}\|\nabla \eta^{j,n}\|^{2} + {\frac{C \, C_{b^*}^2 \, C_{stab}}{\epsilon \, \nu_{max}}}\|\nabla \overline{u_{R}^{j,n}}\| \|\nabla \eta^{j,n+1}\|^{2}
	\\ & {\frac{C \, C_{b^*}^2}{\epsilon \nu_{max}}\|\nabla \eta^{j,n}\|^{2}}+ {\frac{C\,C_{b^*}^2 \, \Delta t }{\epsilon \, \nu_{max}}} \|\nabla(\overline{u_{R}^{j, n}-<u_{R}
		>^{n}})\|^{2}
	\\ &+ {\frac{C}{d \, \epsilon \, \nu_{max}}} \|p^{j,n+1} - q_{h}^{n+1} \|^{2} + {\frac{C \, \Delta t} {\epsilon \, \nu_{max} } }
	\\ &+  {\frac{C\,C_{b^*}^2 \, \Delta t }{\epsilon \, \nu_{max}}} +  {\frac{C\,C_{b^*}^2 \,}{\epsilon \, \nu_{max}}} \|u_{R}^{j,n} - \overline{u_{R}^{j,n}}\| \|\nabla (u_{R}^{j,n} - \overline{u_{R}^{j,n}})\|\|\nabla u^{j,n+1}\|^{2} \bigg \}.
	\end{aligned}
	\end{equation}
	Now using assumption \ref{assumption1}, lemma \ref{lm:ROMstability}, and the stability result from theorem \ref{stab:theorem}, i.e. {$\frac{\epsilon \, \nu_{max} \Delta t}{4} \sum_{n=0}^{N-1} \|\nabla u_{R}^{j,n} \|^{2} \leq C_{stab}$}, we have
	\begin{equation}
	\begin{aligned}
	  & {\frac{\Delta t\, C \, C_{b^*}^2 \, C_{stab}}{\epsilon \, \nu_{max}}} \sum_{n=0}^{N-1}\|\overline{\nabla u_{R}^{j,n}}\|\| \nabla \eta^{j,n+1}\|^{2} 
	\\ & \leq {\frac{\Delta t\, C \, C_{b^*}^2 \, C_{stab}}{\epsilon \, \nu_{max}}} \left(\sum_{n=0}^{N-1} \frac{1}{2} + \sum_{n=0}^{N-1}\frac{\|\nabla u_{R}^{j,n}\|^{2}}{2}  \right) \times
	\\ &\biggr((C+h^2 \|\hspace{-1pt}|  {\mathbb S}_R \|\hspace{-1pt}| _2 ) h^{2s} + (C+\|\hspace{-1pt}|  {\mathbb S}_R \|\hspace{-1pt}| _2 )\Delta t^4  
	+ J_S\sum_{i=R+1}^{J_S(N_S+1)} \|  \nabla \varphi_i\| ^2\lambda_i  \biggr).
	\end{aligned}
	\end{equation}
	Rearranging the first term
	\begin{equation}
	\begin{aligned}
	&{{\frac{\Delta t\, C \, C_{b^*}^2 \, C_{stab}}{\epsilon \, \nu_{max}}} \left(\sum_{n=0}^{N-1} \frac{1}{2} + \sum_{n=0}^{N-1}\frac{\|\nabla u_{R}^{j,n}\|^{2}}{2}  \right)} 
	\\ &= {{\frac{C \, C_{b^*}^2 \, C_{stab}}{2 \, \epsilon \, \nu_{max}}} + {\frac{2 \,C \, C_{b^*}^2 \, C_{stab}}{\epsilon^{2} \, \nu^{2}_{max}}} \frac{\epsilon \nu_{max} \Delta t}{4} \left(\sum_{n=0}^{N-1}\|\nabla u_{R}^{j,n}\|^{2}  \right)}.
	\end{aligned}
	\end{equation}
	It then follows that
	\begin{equation}
	\begin{aligned}
	&{\frac{\Delta t\, C \, C_{b^*}^2 \, C_{stab}}{\epsilon \, \nu_{max}} \sum_{n=0}^{N-1}\|\overline{\nabla u_{R}^{j,n}}\|\| \nabla \eta^{j,n+1}\|^{2}} 
	\\ & {\leq \left(\frac{C \, C_{b^*}^2 \, C_{stab}}{2 \, \epsilon \, \nu_{max}} + \frac{2 \,C \, C_{b^*}^2 \, C^{2}_{stab}}{\epsilon^{2} \, \nu^{2}_{max}} \right) \times}
	\\ &{\biggr((C+h^2 \|\hspace{-1pt}|  {\mathbb S}_R \|\hspace{-1pt}| _2 ) h^{2s} + (C+\|\hspace{-1pt}|  {\mathbb S}_R \|\hspace{-1pt}| _2 )\Delta t^4  
	+ J_S\sum_{i=R+1}^{J_S(N_S+1)} \|  \nabla \varphi_i\| ^2\lambda_i  \biggr)}.
	\end{aligned}
	\end{equation}
	{Next using lemma \ref{lm:ROMfiltering} and Assumptions \ref{assumption:reg} and \ref{assumption:leray}}
	\begin{equation}
	\label{eqn:theorem-convergence-1}
	\begin{aligned}
	&  {\frac{C\,C_{b^*}^2 \, \Delta t}{\epsilon \, \nu_{max}}} \sum_{n=0}^{N-1} \|u_{R}^{j,n} - \overline{u_{R}^{j,n}}\| \|\nabla (u_{R}^{j,n} - \overline{u_{R}^{j,n}})\|\|\nabla u^{j,n+1}\|^{2} 
	\\ &\leq {\frac{C\,C_{b^*}^2 \, \Delta t}{\epsilon \, \nu_{max}}} \sum_{n=0}^{N-1} \|\nabla u^{j,n+1}\|^{2} \frac{1}{\delta} \biggr[C\biggr(C\left( h^{2s+2} + \Delta t^4  \right)+ J_S\sum_{i=R+1}^{J_S(N_S+1)} \lambda_i  \biggr) 
	\\& + C\delta^{2}  \biggr((C+h^2 \|\hspace{-1pt}|  {\mathbb S}_R \|\hspace{-1pt}| _2 ) h^{2s} + (C+\|\hspace{-1pt}|  {\mathbb S}_R \|\hspace{-1pt}| _2 )\Delta t^4  
	+ J_S\sum_{i=R+1}^{J_S(N_S+1)} \|  \nabla \varphi_i\| ^2\lambda_i \biggr) 
	\\& + C\delta^{4}\|\Delta u_{R}^{j,n} \|^{2}\biggr]
	\\ &\leq {\frac{C\,C_{b^*}^2}{\epsilon \, \nu_{max} \, \delta}} \biggr[C\biggr(C\left( h^{2s+2} + \Delta t^4  \right)+ J_S\sum_{i=R+1}^{J_S(N_S+1)} \lambda_i  \biggr) 
	\\& + C\delta^{2}  \biggr((C+h^2 \|\hspace{-1pt}|  {\mathbb S}_R \|\hspace{-1pt}| _2 ) h^{2s} + (C+\|\hspace{-1pt}|  {\mathbb S}_R \|\hspace{-1pt}| _2 )\Delta t^4  
	+ J_S\sum_{i=R+1}^{J_S(N_S+1)} \|  \nabla \varphi_i\| ^2\lambda_i \biggr) 
	\\& + C\delta^{4}\|\Delta u_{R}^{j,n} \|^{2}\biggr].
	\end{aligned}
	\end{equation}
	Next using theorem \ref{stab:theorem} we have
	\begin{equation}
	\begin{aligned}
	\frac{\Delta t \, C\, C_{b^{*}}^{2}}{\epsilon \, \nu_{max}} \|\nabla(\overline{u_{R}^{j, n}-<u_{R} >^{n}})\|^{2} &= 
	\frac{C \||\mathbb{S}_{R}|\|_{2}^{-\frac{1}{2}} }{\epsilon}\frac{C_{b^{*}}^{2} \Delta t }{\nu_{max}}\||\mathbb{S}_{R}|\|_{2}^{\frac{1}{2}} \|\nabla(\overline{u_{R}^{j, n}-<u_{R} >^{n}})\|^{2}
	\\ &\leq C \||\mathbb{S}_{R}|\|_{2}^{-\frac{1}{2}}.
	\end{aligned}
	\end{equation}
	Therefore we can bound the quantities 
	\begin{equation}
	\begin{aligned}
	&\frac{\Delta t \, C\, C_{b^{*}}^{2}}{\epsilon \, \nu_{max}} \|\nabla(\overline{u_{R}^{j, n}-<u_{R} >^{n}})\|^{2} \|\eta_{j}^{n+1}\|^{2} \leq  C \||\mathbb{S}_{R}|\|_{2}^{-\frac{1}{2}} \|\eta_{j}^{n+1}\|^{2} \\
	&\frac{\Delta t \, C\, C_{b^{*}}^{2}}{\epsilon \, \nu_{max}} \|\nabla(\overline{u_{R}^{j, n}-<u_{R} >^{n}})\|^{2} \|\eta_{j}^{n}\|^{2} \leq  C \||\mathbb{S}_{R}|\|_{2}^{-\frac{1}{2}} \|\eta_{j}^{n}\|^{2} \\
	&\frac{\Delta t^{2} \, C\, C_{b^{*}}^{2}}{\epsilon \, \nu_{max}} \|\nabla(\overline{u_{R}^{j, n}-<u_{R} >^{n}})\|^{2}  \leq  C \Delta t \||\mathbb{S}_{R}|\|_{2}^{-\frac{1}{2}}.
	\end{aligned}
	\end{equation}
	%\begin{equation}
	%\frac{\Delta t \, C\, C_{b^{*}}^{2}}{\epsilon \, \nu_{max}} \|\nabla(\overline{u_{R}^{j, n}-<u_{R} >^{n}})\|^{2} \|\eta_{j}^{n+1}\|^{2} = \frac{\||\mathbb{S}_{R}|\|_{2}^{\frac{1}{2}}}{\||\mathbb{S}_{R}|\|_{2}^{\frac{1}{2}}}\frac{\Delta t \, C\, C_{b^{*}}^{2}}{\epsilon \, \nu_{max}} \|\nabla(\overline{u_{R}^{j, n}-<u_{R} >^{n}})\|^{2} \|\eta_{j}^{n+1}\|^{2}
	%\end{equation}
	Now combining everything, absorbing constants, invoking the discrete Gronwall's inequality, using Assumption \ref{assumption1} and the stability estimate \eqref{stab:assumption} \eqref{comb3} becomes
	
	\begin{equation}
	\begin{aligned}
	\label{ineq:final}
	&\frac{1}{2}\|\xi_{R}^{j,N}\|^{2} + \frac{\nu_{max}}{2}\| \nabla \xi_{R}^{j,N}\|^{2} + {\frac{\epsilon}{4} \nu_{max} \| \nabla \xi_{R}^{j,N}\|^{2}} + C\nu_{max}\Delta t \sum_{n=0}^{N-1} \| \nabla \xi_{R}^{j,n+1} \|^{2}
	\\ &\leq 
	\exp\left(\frac{C_{b^{*}}^{4} \, C^{4} \,T}
	{\epsilon^{3}\,\nu_{max}^{3}}\right) 
	\biggr[\biggr(\frac{ C \, \nu_{max} \Delta t}{\epsilon}
	+ \frac{C \Delta t}{\epsilon} \frac{|\nu_j - \nu_{max}|^{2}}{\nu_{max}} 
	+ 2C \||\mathbb{S}_{R}|\|_{2}^{-\frac{1}{2}}
	+ \frac{C \, C_{b^*}^2 \, \Delta t}{\epsilon \, \nu_{max}}
	\\ &+ \frac{C \, C_{b^*}^2 \, C_{stab}}{2 \, \epsilon \, \nu_{max}} 
	+ \frac{2 \,C \, C_{b^*}^2 \, C^{2}_{stab}}{\epsilon^{2} \, \nu^{2}_{max}}
	+ \frac{C \, C_{b^*}^2 \, \delta}{ \epsilon \, \nu_{max}}  \biggr) \times 
	\\& \biggr((C+h^2 \|\hspace{-1pt}|  {\mathbb S}_R \|\hspace{-1pt}| _2 ) h^{2s} + (C+\|\hspace{-1pt}|  {\mathbb S}_R \|\hspace{-1pt}| _2 )\Delta t^4  
	+ J_S\sum_{i=R+1}^{J_S(N_S+1)} \|  \nabla \varphi_i\| ^2\lambda_i \biggr) 
	\\& + \frac{C\,C_{b^*}^2}{\delta \,\epsilon \, \nu_{max}} \biggr(C\left( h^{2s+2} + \Delta t^4  \right)+ J_S\sum_{i=R+1}^{J_S(N_S+1)} \lambda_i  \biggr)
	\\& + \frac{C \, \Delta t^{2}}{\epsilon}\frac{|\nu_j - \nu_{max}|^{2}}{\nu_{max}} 
	+ C \Delta t \||\mathbb{S}_{R}|\|_{2}^{-\frac{1}{2}}
	+ \frac{C h^{2s}} {d \, \epsilon \, \nu_{max}}\| |p^{j} | \|^{2}_{2,s} 
	+ \frac{C \Delta t^{2}}{\epsilon \, \nu_{max}} 
	\\& 
	+ \frac{C\,C_{b^*}^2 \, \Delta t^{2}}{\epsilon \, \nu_{max}} 
	+ \frac{C\,C_{b^*}^2 \, \delta^3}{\epsilon \, \nu_{max} \, } \biggr].
	\end{aligned}
	\end{equation}	
%\biggr( \frac{CN_{S}\Delta t}{\nu_{max}} + CN_{S}\Delta t\nu_{max}  + C N_{S}\Delta t\frac{|\nu_{j} - \nu_{max}|^{2}}{\nu_{max}} + C N_{S} \|\hspace{-1pt}|  {\mathbb S}_R \|\hspace{-1pt}|_{2}^{^{-\frac{1}{2}}} + C \delta \biggr) \cdot
%	\\ & \biggr( \inf_{{j\in\{1,\ldots,J_S\}}}\frac{2}{N_S+1} \sum_{m=0}^{N_S} \left(\|  \nabla (u^{m}-u_S^{j,m})\| ^2+\|\hspace{-1pt}|  {\mathbb S}_R \|\hspace{-1pt}| _2 \|   u^{m}-u_S^{j,m}\| ^2  \right)
%	\\& + (C+h^2 \|\hspace{-1pt}|  {\mathbb S}_R \|\hspace{-1pt}| _2 ) h^{2s} + (C+\|\hspace{-1pt}|  {\mathbb S}_R \|\hspace{-1pt}| _2 )\Delta t^4  
%	+ 2J_S\sum_{i=R+1}^{J_S(N_S+1)} \|  \nabla \varphi_i\| ^2\lambda_i \biggr) 
%	\\ & +  \frac{1}{\delta} \biggr( \inf_{j\in\{1,\ldots,J_S\}}\frac{2}{N_S+1} \sum_{m=0}^{N_S} \|  u^{m}-u^{j,m}_S\| ^2 +C\left( h^{2s+2} + \Delta t^4  \right)+ 2J_S\sum_{i=R+1}^{J_S(N_S+1)} \lambda_i  \biggr) 
%	\\ & + C\delta ^{3} + C\Delta t^{2}\frac{|\nu_{j} - \nu_{max}|^{2}}{\nu_{max}}\| |\nabla u_{j,t}|\|_{2,0}^{2} + \frac{C h^{2s}} {\nu_{max}}\| |p^{j} | \|^{2}_{2,s} + \frac{C\Delta t^{2}}{\nu_{max}}\| |u_{j,tt}| \|^{2}_{2,0} 
%	\\ & +  \frac{C \Delta t^{2}}{\nu_{max}}\| |\nabla u_{j,t}|\|_{2,0}^{2} + C \Delta t\| \|\hspace{-1pt}|  {\mathbb S}_R \|\hspace{-1pt}|_{2}^{^{-\frac{1}{2}}} \||\nabla u_{j,t}|\|_{2,0}^{2}.
%	\end{aligned}
%	\end{equation}

	By the triangle inequality we have $\|e^{j,n}\|^{2} \leq 2(\|\xi_{R}^{j,n} \|^{2} + \| \eta^{j,n}\|^{2})$ from which it follows
	
	\begin{equation}
	\begin{aligned}
	&\frac{1}{2}\|e^{j,N}\|^{2} + \frac{\nu_{max}}{2}\| \nabla e^{j,N}\|^{2} +  {\frac{\epsilon}{4}\nu_{max} \| \nabla e^{j,N}\|^{2}}  + C\nu_{max}\Delta t \sum_{n=0}^{N-1} \| e^{j,n+1} \|^{2} 
	\\ & \leq \|\eta^{j,N}\|^{2} + \nu_{max}\| \nabla \eta^{j,N}\|^{2} + {\frac{\epsilon}{2} \nu_{max} \| \nabla \eta^{j,N}\|^{2}}  + C\nu_{max}\Delta t \sum_{n=0}^{N-1} \| \eta^{j,n+1} \|^{2} 
	\\ & + \|\xi_{R}^{j,N}\|^{2} + \nu_{max}\| \nabla \xi_{R}^{j,N}\|^{2} + {\frac{\epsilon}{2} \nu_{max} \| \nabla \xi_{R}^{j,N}\|^{2}} + C\nu_{max}\Delta t \sum_{n=0}^{N-1} \| \xi_{R}^{j,n+1} \|^{2}. 	
	\end{aligned}
	\end{equation}
	Now applying inequality \eqref{ineq:final} and Assumption \ref{assumption1} the result follows.
\end{proof}
\end{theorem}

\section{Numerical Experiments}\label{numex}
In this section we provide numerical experiments for the Leray ensemble-POD algorithm \eqref{En-Leray-POD-Weak} demonstrating the efficacy of this approach. All computations will be done using the FEniCS software suite \cite{LNW12} and all meshes generated via the built in meshing package \textbf{mshr}.

%\subsection{Differential Filter Implementation}
%In order to implement the differential filter there are two commonly used approaches. (\blue{input information about FE version}), in this work we will use the ROM approach. Using this approach at each time step we need to solve the $r\times r$ reduced basis Helmholtz system
%\begin{equation}
%(M_{R} + \delta^{2}S_{R})\mathcal{F}(\vec{a}) = M_{R}\vec{a}.
%\end{equation}
%Since the system is only dependent on the POD basis it only needs to be assembled and factorized once. It will then reduce at each time step to an $O(r^{2})$ solve. In Figure \ref{filtered_average} we give an example of the impact the filter will have on the average solution at time $T = 1.01$.

\subsection{Problem Setting}
\begin{figure}[h!]
	\centering
	\includegraphics[width = 12cm]{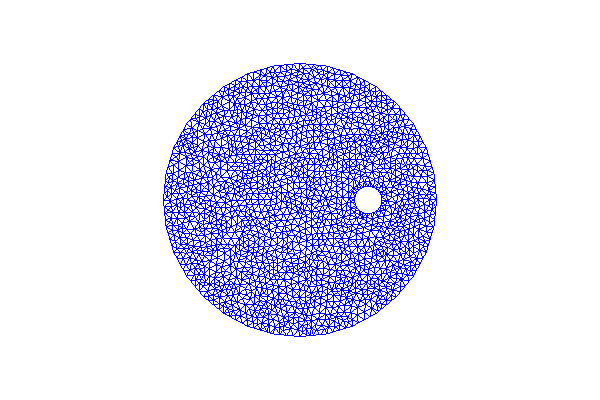} 
	\caption{Mesh for flow between offset circles resulting in 14,590 total degrees of freedom for the Taylor-Hood element pair.}
	\label{meshCinC}
\end{figure}

For the numerical experiments we consider the two-dimensional flow between offset cylinders used in \cite{GJS17}. The domain is a disk with a smaller off-center disc inside. Let $r_{1}=1$, $r_{2}=0.1$, $c_{1}=1/2$, and $c_{2}=0$; then, the domain is given by
\[
\Omega=\{(x,y):x^{2}+y^{2}\leq r_{1}^{2} \text{ and } (x-c_{1})^{2}%
+(y-c_{2})^{2}\geq r_{2}^{2}\}.
\]
The mesh utilized contains 14,590 degrees of freedom and is given in figure \ref{meshCinC}. We discretize in space via the $P^2$-$P^1$ Taylor-Hood element pair. The no-slip, no-penetration boundary conditions are imposed on both cylinders. In our test problems the flow will be driven by the counterclockwise rotational body force
\[
f(x,y,t)=\big(-4y(1-x^{2}-y^{2})\,,\,4x(1-x^{2}-y^{2})\big)^{T}.
\]
This flow displays interesting structures which interact with the inner circle. Specifically the flow rotates about the origin and interacts with the immersed cylinder forming a Von K\'arm\'an vortex street.

\subsection{Numerical Results}
\begin{figure}[h!]
	\centering
	\includegraphics[width=10cm]{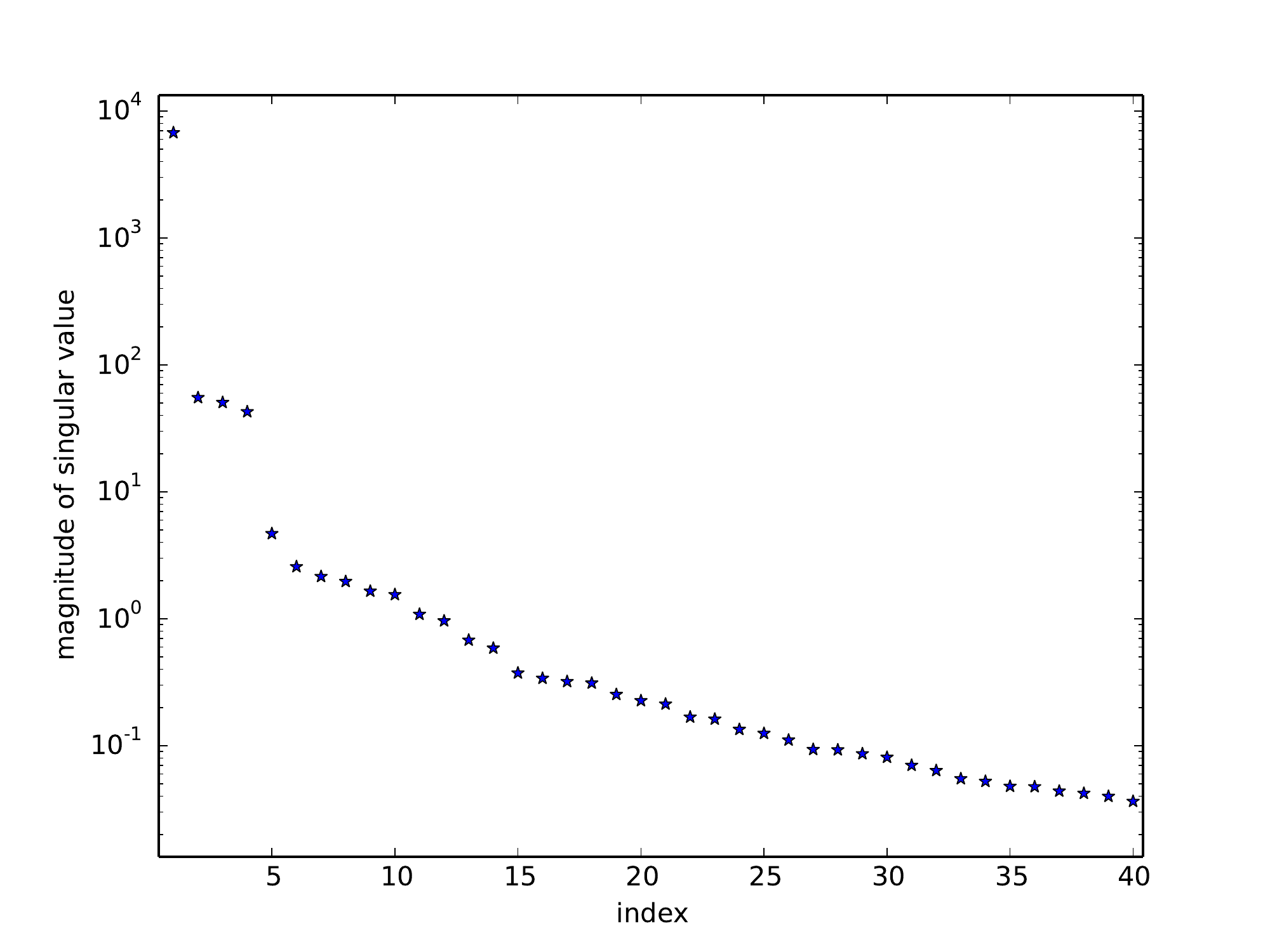} 
	\caption{The 40 largest eigenvalues for the snapshot matrix.}
	\label{eigvals_ex1}
\end{figure}

In this experiment we demonstrate the improved accuracy and stability of the Leray ensemble-POD algorithm. In order to generate the POD basis we use two different viscosities $\nu_1 = .0016$ and $\nu_2 = .002$. The initial conditions will be generated by solving a steady Stokes problem using the previously defined counterclockwise rotational body force. We run a finite element code utilizing a linearly implicit backwards Euler method for each viscosity from $t_0 = 0$ to $T= 6$ with fixed time step $\Delta t = .01$. At time $T=3.0$  we begin taking snapshots every $.04$ seconds. In Figure \ref{eigvals_ex1} we show the decay of the singular values for the snapshot matrix.

To illustrate the accuracy of the Leray ensemble-POD algorithm we compare it against the ensemble-POD  algorithm using the same viscosities from the offline stage. The computations are carried out over the time interval $t_0 = 3$ to $T= 6$ with fixed time step $\Delta t = .01$ and $r = 10$ reduced basis functions. The initial condition at $T = 3$ is the $L^{2}$ projection of the FE solution at $T = 3.0$ into the POD space The filtering length for the Leray ensemble-POD algorithm is taken to be $\delta = .025$. The filtering length is selected as the value of $\delta$ which allows the average kinetic energy of Leray ensemble-POD to most closely match the average kinetic energy of the benchmark solution. We purposefully utilize a small number of basis functions to demonstrate the situation where the ROM does not allow for all spatial scales to be resolved. To determine the accuracy of our methods the average of the solutions from the implicit backwards Euler method for $\nu_1$ and $\nu_2$ will be used as a benchmark. 

\begin{figure}[h!]
\centering
\includegraphics[width=12.cm]{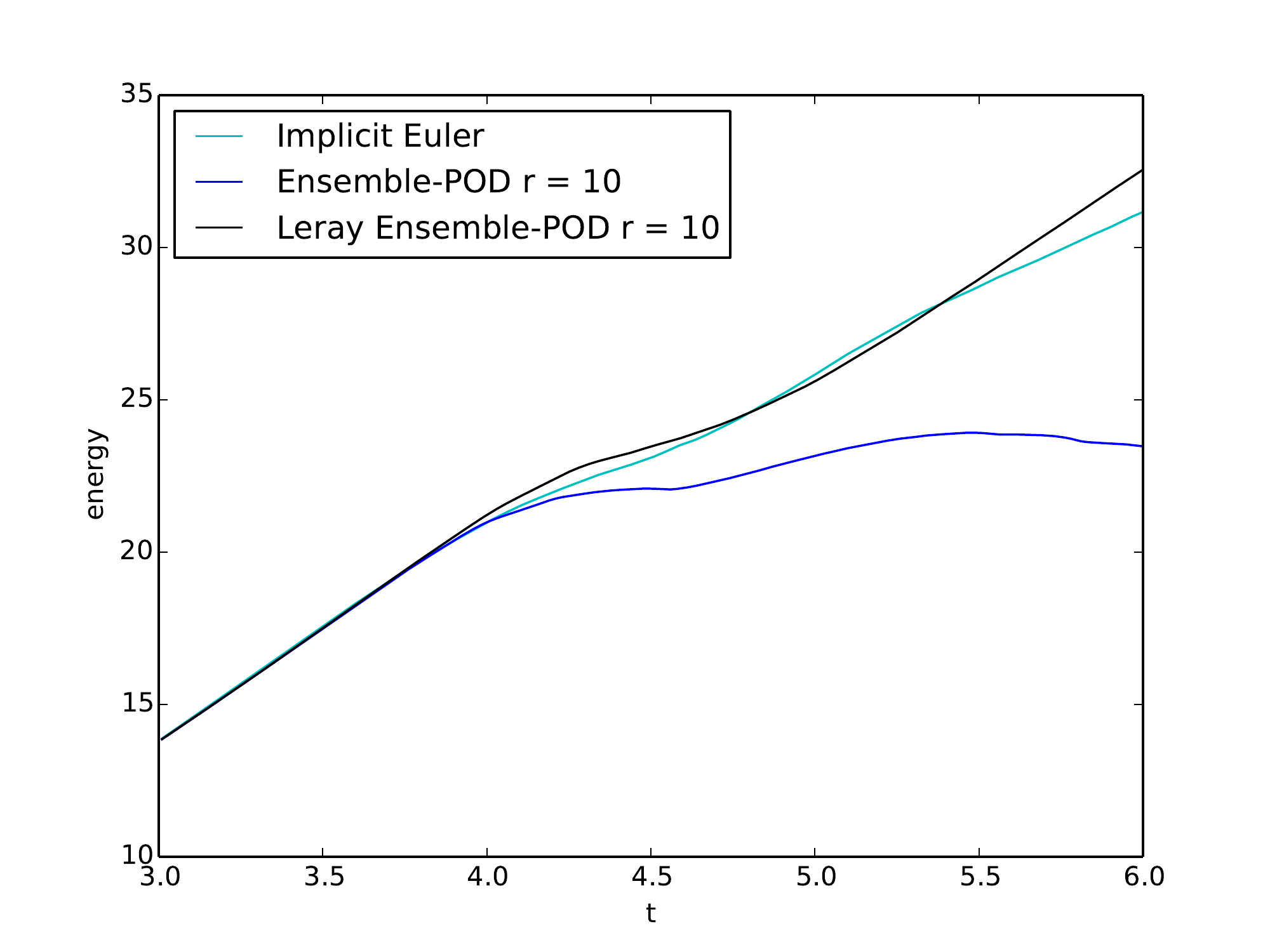}
\caption{For $3 \leq t \leq 6$, the average energy of the Leray ensemble-POD, ensemble-POD, and Implicit Euler method.}
\label{KE_fig}
\end{figure}

\begin{figure}[h!]
\centering
\includegraphics[width=12.cm]{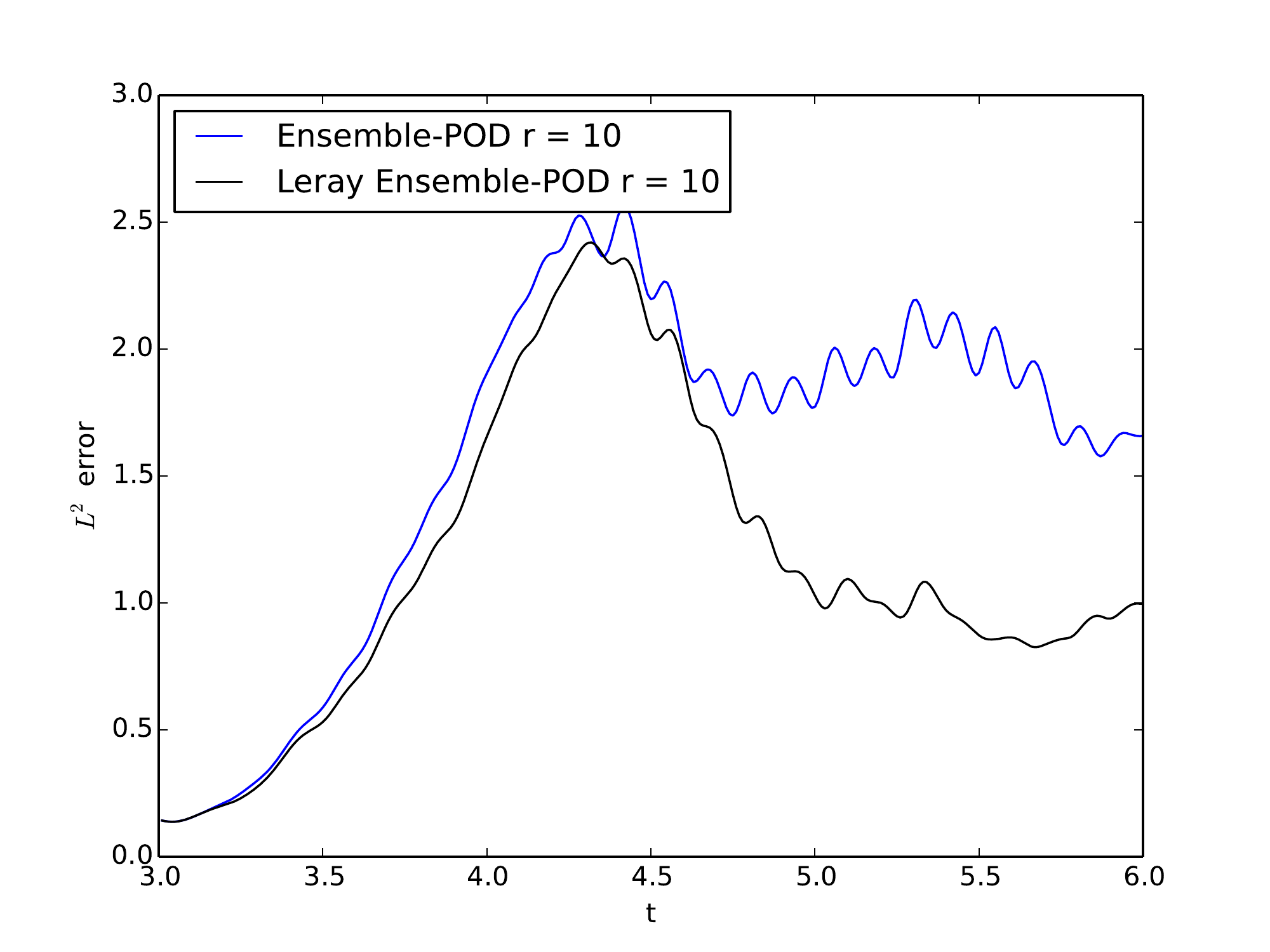}
\caption{For $3 \leq t \leq 6$, the $L^{2}$ error evolution of the Leray ensemble-POD and ensemble-POD algorithms with $r=10$ basis functions.}
\label{error_fig}
\end{figure}

\begin{figure}[h!]
\centering
\includegraphics[width=6.cm]{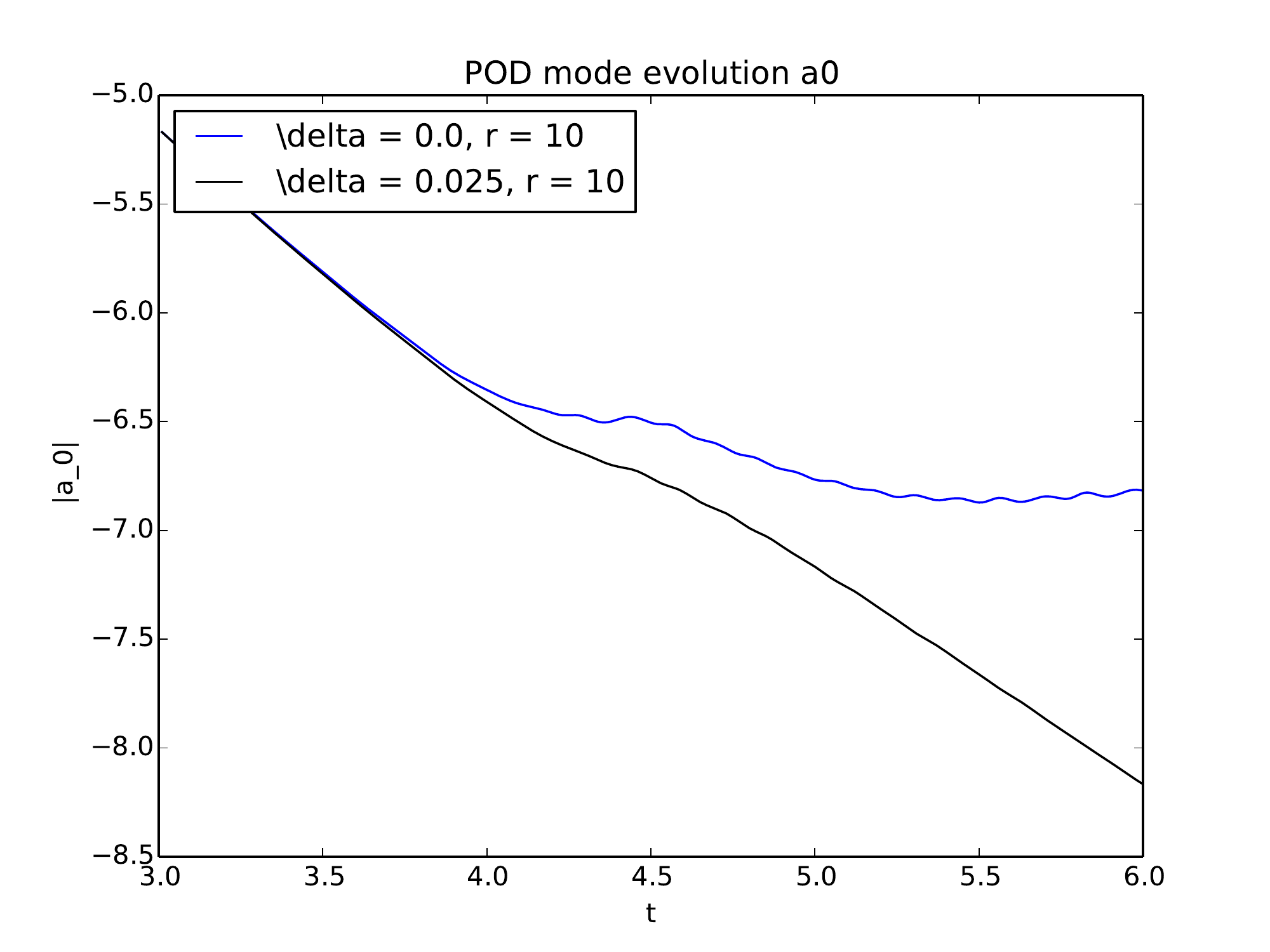}
\includegraphics[width=6.cm]{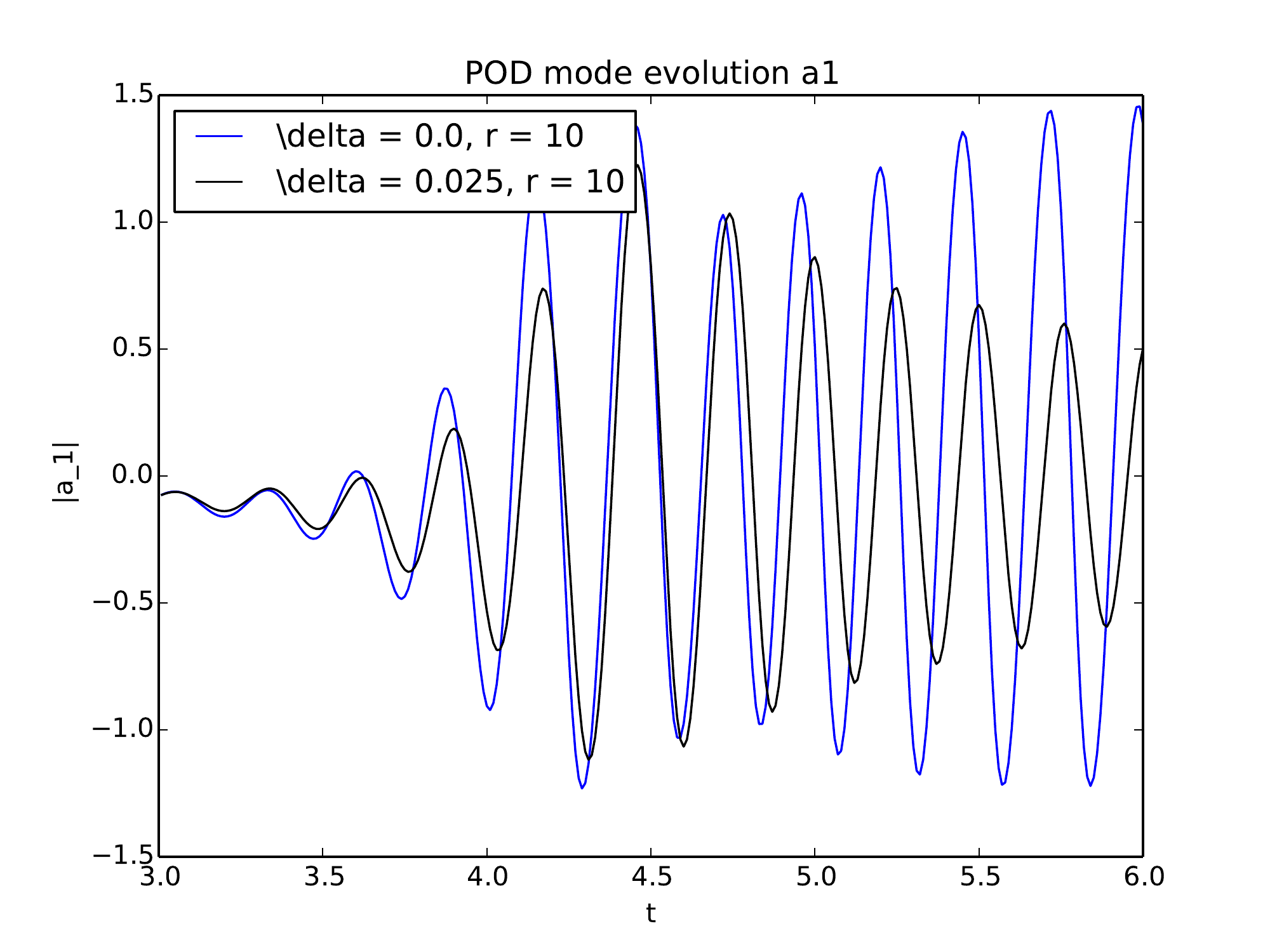} 
\includegraphics[width=6.cm]{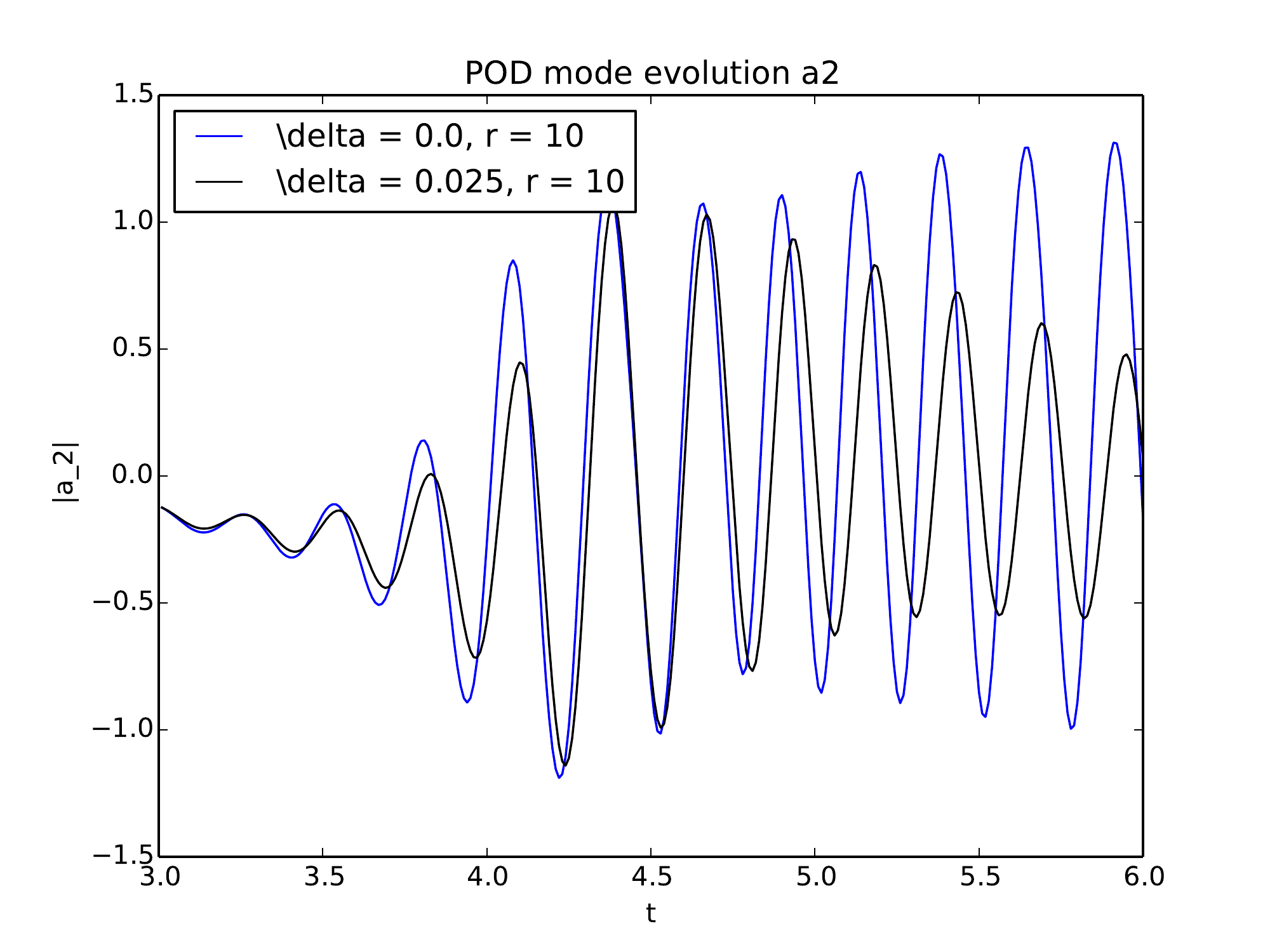} 
\includegraphics[width=6.cm]{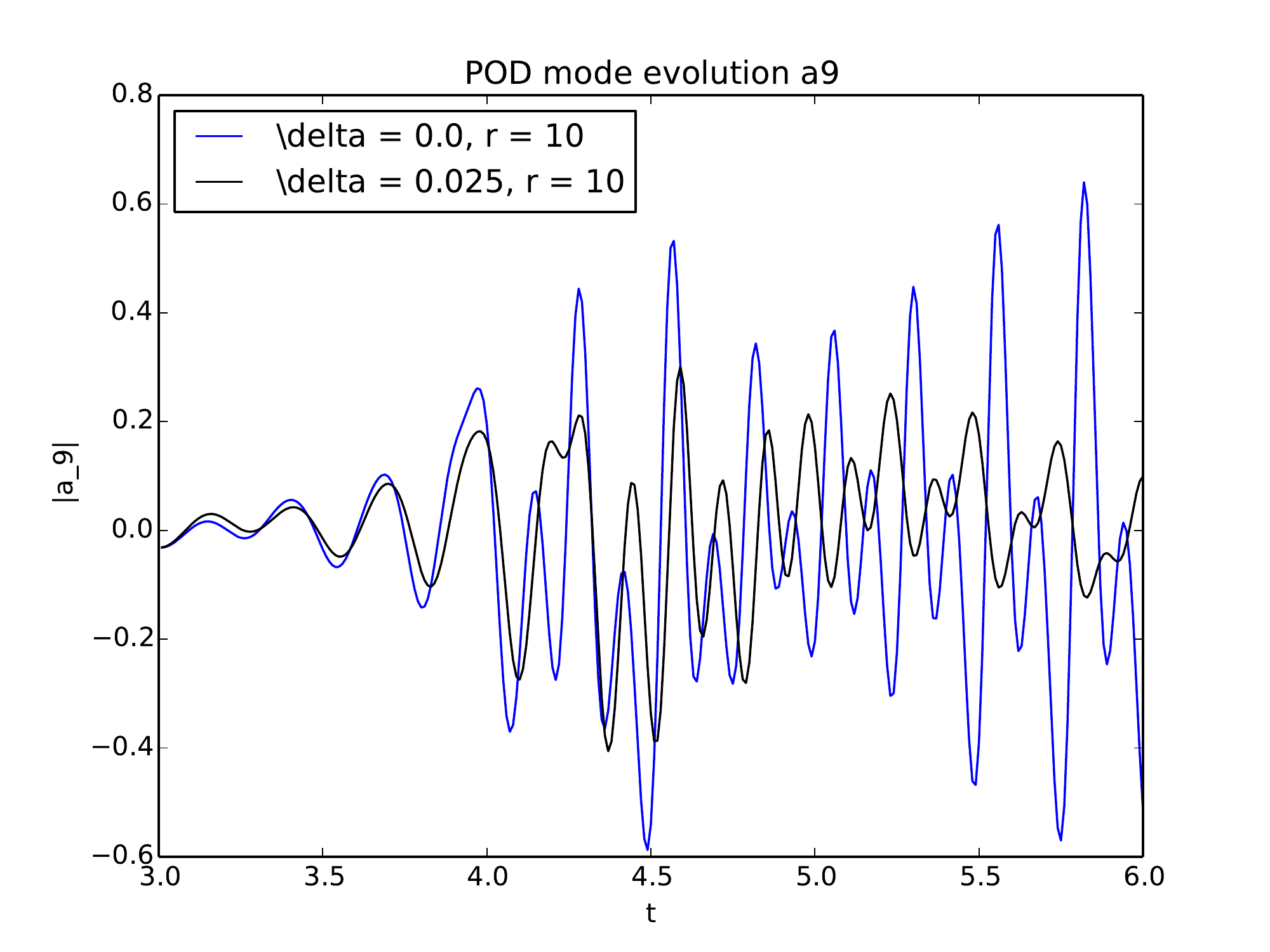} 
\caption{The time evolution of the average POD modes for ensemble-POD versus Leray ensemble-POD.}
\label{mode_evolution}
\end{figure}

In Figure \ref{KE_fig} we compare the average kinetic energy evolution of the Leray ensemble-POD and ensemble-POD against our benchmark solution. It can be seen that the ensemble-POD fails to match the kinetic energy of the benchmark solution, while the Leray ensemble-POD approximates it reasonably well. In Figure \ref{error_fig} we compare the evolution of the error in the $L^{2}$ norm of Leray ensemble-POD and ensemble-POD algorithms. The Leray ensemble-POD has a significantly smaller error than the ensemble-POD algorithm. In Figure \ref{mode_evolution} we plot the average POD mode evolution for ensemble-POD versus Leray ensemble-POD. We see that the oscillations in the POD modes are damped for Leray ensemble-POD.

\section{Conclusions}
In this work, a Leray regularized ensemble-POD method is developed for the incompressible Navier-Stokes equations with perturbations in the forcing function, initial conditions, and viscosities. The proposed algorithm is first ensemble-POD approach designed to work for higher Reynolds number flows. The stability and convergence of the finite element discretization of the Leray ensemble-POD model are proven.  In the numerical simulation of two-dimensional flow past two offset cylinders, it is shown that the Leray ensemble-POD model is significantly more accurate than the standard ensemble-POD model.
\bibliography{WorksCited}

\begin{thebibliography}{10}

\bibitem{balajewicz2013low}
M.~J. Balajewicz, E.~H. Dowell, and B.~R. Noack.
\newblock Low-dimensional modelling of high-{R}eynolds-number shear flows
  incorporating constraints from the {N}avier--{S}tokes equation.
\newblock {\em J. Fluid Mech.}, 729:285--308, 2013.

\bibitem{balajewicz2016minimal}
M.~J. Balajewicz, I.~Tezaur, and E.~H. Dowell.
\newblock {Minimal subspace rotation on the Stiefel manifold for stabilization
  and enhancement of projection-based reduced order models for the compressible
  Navier--Stokes equations}.
\newblock {\em J. Comput. Phys.}, 321:224--241, 2016.

\bibitem{NME:NME4772}
Francesco Ballarin, Andrea Manzoni, Alfio Quarteroni, and Gianluigi Rozza.
\newblock Supremizer stabilization of pod–galerkin approximation of
  parametrized steady incompressible navier–stokes equations.
\newblock {\em International Journal for Numerical Methods in Engineering},
  102(5):1136--1161, 2015.

\bibitem{benosman2016robust}
M.~Benosman, J.~Borggaard, and B.~Kramer.
\newblock Robust reduced-order model stabilization for partial differential
  equations based on {L}yapunov theory and extremum seeking with application to
  the {3D} {B}oussinesq equations.
\newblock {\em arXiv preprint arXiv:1604.04586}, 2016.

\bibitem{carlberg2013gnat}
K.~Carlberg, C.~Farhat, J.~Cortial, and D.~Amsallem.
\newblock The {GNAT} method for nonlinear model reduction: effective
  implementation and application to computational fluid dynamics and turbulent
  flows.
\newblock {\em J. Comput. Phys.}, 2013.

\bibitem{2017arXiv170900243C}
T.~{Chac{\'o}n Rebollo}, E.~{Delgado {\'A}vila}, M.~{G{\'o}mez M{\'a}rmol},
  F.~{Ballarin}, and G.~{Rozza}.
\newblock {On a certified Smagorinsky reduced basis turbulence model}.
\newblock {\em ArXiv e-prints}, September 2017.

\bibitem{FOP95}
YT~Feng, DRJ Owen, and D~Peri{\'c}.
\newblock A block conjugate gradient method applied to linear systems with
  multiple right-hand sides.
\newblock {\em Computer methods in applied mechanics and engineering},
  127(1):203--215, 1995.

\bibitem{2017arXiv171003569F}
L.~{Fick}, Y.~{Maday}, A.~T {Patera}, and T.~{Taddei}.
\newblock {A Reduced Basis Technique for Long-Time Unsteady Turbulent Flows}.
\newblock {\em ArXiv e-prints}, October 2017.

\bibitem{FM97}
Roland~W Freund and Manish Malhotra.
\newblock A block qmr algorithm for non-hermitian linear systems with multiple
  right-hand sides.
\newblock {\em Linear Algebra and its Applications}, 254(1):119--157, 1997.

\bibitem{GS96}
E~Gallopulos and V~Simoncini.
\newblock Convergence of block gmres and matrix polynomials.
\newblock {\em Linear Algebra Appl}, 247:97--119, 1996.

\bibitem{germano1986differential}
M.~Germano.
\newblock Differential filters of elliptic type.
\newblock {\em Phys. Fluids}, 29(6):1757--1758, 1986.

\bibitem{geurts2003regularization}
B.~J. Geurts and D.~D. Holm.
\newblock Regularization modeling for large-eddy simulation.
\newblock {\em Phys. Fluids}, 15(1):L13--L16, 2003.

\bibitem{giere2015supg}
S.~Giere, T.~Iliescu, V.~John, and D.~Wells.
\newblock {SUPG} reduced order models for convection-dominated
  convection-diffusion-reaction equations.
\newblock {\em Comput. Methods Appl. Mech. Engrg.}, 289:454--474, 2015.

\bibitem{GR79}
V.~Girault and P.-A. Raviart.
\newblock {\em Finite element approximation of the {N}avier-{S}tokes
  equations}, volume 749 of {\em Lecture Notes in Mathematics}.
\newblock Springer-Verlag, Berlin, 1979.

\bibitem{2017arXiv170604060G}
M.~{Gunzburger}, N.~{Jiang}, and Z.~{Wang}.
\newblock {A second-order time-stepping scheme for simulating ensembles of
  parameterized flow problems}.
\newblock {\em ArXiv e-prints}, June 2017.

\bibitem{GJW17}
M.~{Gunzburger}, N.~{Jiang}, and Z.~{Wang}.
\newblock {An efficient algorithm for simulating ensembles of parameterized
  flow problems}.
\newblock {\em ArXiv e-prints}, May 2017.

\bibitem{GJS17}
Max Gunzburger, Nan Jiang, and Michael Schneier.
\newblock An ensemble-proper orthogonal decomposition method for the
  nonstationary navier-stokes equations.
\newblock {\em SIAM Journal on Numerical Analysis}, 55(1):286--304, 2017.

\bibitem{GJS16}
Max Gunzburger, Nan Jiang, and Michael Schneier.
\newblock A higher-order ensemble/proper orthogonal decomposition method for
  the nonstationary navier-stokes equations.
\newblock {\em International Journal of Numerical Analysis and Modeling}, to
  appear, 2018.

\bibitem{Max89}
Max~D Gunzburger.
\newblock {\em Finite Element Methods for Viscous Incompressible Flows: A guide
  to theory, practice, and algorithms}.
\newblock Elsevier, 2012.

\bibitem{gunzburger_webster_zhang_2014}
Max~D. Gunzburger, Clayton~G. Webster, and Guannan Zhang.
\newblock Stochastic finite element methods for partial differential equations
  with random input data.
\newblock {\em Acta Numerica}, 23:521–650, 2014.

\bibitem{hesthaven2015certified}
J.~S. Hesthaven, G.~Rozza, and B.~Stamm.
\newblock {\em Certified Reduced Basis Methods for Parametrized Partial
  Differential Equations}.
\newblock Springer, 2015.

\bibitem{iliescu2017regularized}
T.~{Iliescu}, H.~{Liu}, and X.~{Xie}.
\newblock {Regularized Reduced Order Models for a Stochastic Burgers Equation}.
\newblock {\em International Journal of Numerical Analysis and Modeling}, to
  appear, 2018.

\bibitem{IW14}
Traian Iliescu and Zhu Wang.
\newblock Variational multiscale proper orthogonal decomposition: Navier-stokes
  equations.
\newblock {\em Numerical Methods for Partial Differential Equations},
  30(2):641--663, 2014.

\bibitem{J15}
Nan Jiang.
\newblock A higher order ensemble simulation algorithm for fluid flows.
\newblock {\em Journal of Scientific Computing}, pages 1--25, 2014.

\bibitem{J17}
Nan Jiang.
\newblock A second-order ensemble method based on a blended backward
  differentiation formula timestepping scheme for time-dependent navier stokes
  equations.
\newblock {\em Numerical Methods for Partial Differential Equations},
  33(1):34--61, 2017.

\bibitem{JL14}
Nan Jiang and William Layton.
\newblock An algorithm for fast calculation of flow ensembles.
\newblock {\em International Journal for Uncertainty Quantification}, 4(4),
  2014.

\bibitem{john2016divergence}
V.~John, A.~Linke, C.~Merdon, M.~Neilan, and L.~G. Rebholz.
\newblock On the divergence constraint in mixed finite element methods for
  incompressible flows.
\newblock {\em SIAM Rev.}, 2016.

\bibitem{kalashnikova2010stability}
I.~Kalashnikova and M.~F. Barone.
\newblock On the stability and convergence of a {G}alerkin reduced order model
  {(ROM)} of compressible flow with solid wall and far-field boundary
  treatment.
\newblock {\em Int. J. Num. Meth. Eng.}, 83(10):1345--1375, 2010.

\bibitem{KLV01}
Karl Kunisch and Stefan Volkwein.
\newblock Galerkin proper orthogonal decomposition methods for parabolic
  problems.
\newblock {\em Numerische mathematik}, 90(1):117--148, 2001.

\bibitem{layton2008numerical}
W.~Layton, C.~C. Manica, M.~Neda, and L.~G. Rebholz.
\newblock Numerical analysis and computational testing of a high accuracy
  {L}eray-deconvolution model of turbulence.
\newblock {\em Num. Meth. P.D.E.s}, 24(2):555--582, 2008.

\bibitem{layton2012approximate}
W.~J. Layton and L.~G. Rebholz.
\newblock {\em Approximate Deconvolution Models of Turbulence: Analysis,
  Phenomenology and Numerical Analysis}, volume 2042.
\newblock Springer Berlin Heidelberg, 2012.

\bibitem{Layton08}
William Layton.
\newblock {\em Introduction to the numerical analysis of incompressible viscous
  flows}, volume~6.
\newblock Siam, 2008.

\bibitem{leray1934sur}
J.~Leray.
\newblock Sur le mouvement d`un fluide visqueux emplissant l'espace.
\newblock {\em Acta Math.}, 63:193--248, 1934.

\bibitem{LNW12}
Anders Logg, Kent-Andre Mardal, and Garth Wells.
\newblock {\em Automated solution of differential equations by the finite
  element method: The FEniCS book}, volume~84.
\newblock Springer Science \& Business Media, 2012.

\bibitem{LW17}
Y.~Luo and Z.~Wang.
\newblock {An ensemble algorithm for numerical solutions to deterministic and
  random parabolic PDEs}.
\newblock {\em ArXiv e-prints}, October 2017.

\bibitem{MohebujjamanR17}
Muhammad Mohebujjaman and Leo~G. Rebholz.
\newblock An efficient algorithm for computation of {MHD} flow ensembles.
\newblock {\em Comput. Meth. in Appl. Math.}, 17(1):121--137, 2017.

\bibitem{osth2014need}
J.~{\"O}sth, B.~R. Noack, S.~Krajnovi{\'c}, D.~Barros, and J.~Bor{\'e}e.
\newblock On the need for a nonlinear subscale turbulence term in {POD} models
  as exemplified for a high-{R}eynolds-number flow over an {A}hmed body.
\newblock {\em J. Fluid Mech.}, 747:518--544, 2014.

\bibitem{quarteroni2015reduced}
A.~Quarteroni, A.~Manzoni, and F.~Negri.
\newblock {\em Reduced Basis Methods for Partial Differential Equations: An
  Introduction}, volume~92.
\newblock Springer, 2015.

\bibitem{sabetghadam2012alpha}
F.~Sabetghadam and A.~Jafarpour.
\newblock $\alpha$ regularization of the {POD-G}alerkin dynamical systems of
  the {K}uramoto--{S}ivashinsky equation.
\newblock {\em Appl. Math. Comput.}, 218(10):6012--6026, 2012.

\bibitem{AMJ16}
Aziz Takhirov, Monika Neda, and Jiajia Waters.
\newblock Time relaxation algorithm for flow ensembles.
\newblock {\em Numerical Methods for Partial Differential Equations},
  32(3):757--777, 2016.

\bibitem{wang20162d}
Y.~Wang, I.~M. Navon, X.~Wang, and Y.~Cheng.
\newblock {2D Burgers equation with large Reynolds number using POD/DEIM and
  calibration}.
\newblock {\em Int. J. Num. Meth. Fluids}, 82(12):909--931, 2016.

\bibitem{weller2009numerical}
J.~Weller, E.~Lombardi, M.~Bergmann, and A.~Iollo.
\newblock {Numerical methods for low-order modeling of fluid flows based on
  POD}.
\newblock {\em Comput. Meth. Appl. Mech. Eng.}, 200(33-36):2507--2527, 2009.

\bibitem{weller2009robust}
J.~Weller, E.~Lombardi, and A.~Iollo.
\newblock Robust model identification of actuated vortex wakes.
\newblock {\em Phys. D}, 238(4):416--427, 2009.

\bibitem{wells2017evolve}
D.~Wells, Z.~Wang, X.~Xie, and T.~Iliescu.
\newblock An evolve-then-filter regularized reduced order model for
  convection-dominated flows.
\newblock {\em Int. J. Num. Meth. Fluids}, 84:598--–615, 2017.

\bibitem{xie2017data}
X.~Xie, M.~Mohebujjaman, L.~G. Rebholz, and T.~Iliescu.
\newblock Data-driven filtered reduced order modeling.
\newblock {\em arXiv preprint arXiv:1702.06886}, 2017.

\bibitem{xie2017approximate}
X.~Xie, D.~Wells, Z.~Wang, and T.~Iliescu.
\newblock Approximate deconvolution reduced order modeling.
\newblock {\em Comput. Methods Appl. Mech. Engrg.}, 313:512--534, 2017.

\bibitem{XIE201812}
Xuping Xie, David Wells, Zhu Wang, and Traian Iliescu.
\newblock Numerical analysis of the leray reduced order model.
\newblock {\em Journal of Computational and Applied Mathematics},
  328(Supplement C):12 -- 29, 2018.

\end{thebibliography}
\bibliographystyle{plain}

\end{document}